\documentclass[11pt]{amsart}
\usepackage[utf8]{inputenc}
\usepackage{amsmath,amsthm,amsfonts,amssymb,mathtools}
\usepackage{mathrsfs}

\usepackage{kotex}
\usepackage{hyperref}
\usepackage{relsize}

\newcommand{\Be}{\begin{equation}}
\newcommand{\Ee}{\end{equation}}
\newcommand{\Bea}{\begin{eqnarray}}
\newcommand{\Eea}{\end{eqnarray}}
\newcommand{\Bel}{\begin{align}}
\newcommand{\Eel}{\end{align}}
\newcommand{\Beas}{\begin{eqnarray*}}
\newcommand{\Eeas}{\end{eqnarray*}}
\newcommand{\Benu}{\begin{enumerate}}
\newcommand{\Eenu}{\end{enumerate}}
\newcommand{\Bi}{\begin{itemize}}
\newcommand{\Ei}{\end{itemize}}
\newcommand\supp{\operatorname{supp}}

\def\R{{\mathbb R}}
\def\Z{{\mathbb Z}}

\newcommand{\E}{\mathlarger{ \mathcal E}}

\newcommand{\dist}{\mathrm{dist }}

\theoremstyle{plain}
\newtheorem{thm}{Theorem}[section]
\newtheorem{cor}[thm]{Corollary}
\newtheorem{lem}[thm]{Lemma}
\newtheorem{prop}[thm]{Proposition}

\theoremstyle{remark}
  
\theoremstyle{definition}
\newtheorem*{defn}{Definition}

\numberwithin{equation}{section}

\newcommand{\RNum}[1]{\uppercase\expandafter{\romannumeral #1\relax}}

\usepackage{tikz}

\newcommand{\dsig}{{\sigma_{\! t,s}^\theta}}
\newcommand{\dsi}{{\sigma_{\! t,s}^0}}

\newcommand{\ED}{\end{document}}

 \newcommand{\I}{\mathbb I}
\newcommand{\M}{\mathfrak M}
\newcommand{\dil}[1]{\mathrm  D_{\!#1}}
\newcommand{\B}{\mathbb B}

\makeatletter
\@namedef{subjclassname@2020}{\textup{2020} Mathematics Subject Classification}
\makeatother

\title{The elliptic maximal function}
\author[Juyoung Lee]{Juyoung Lee}
\author[Sanghyuk Lee]{Sanghyuk Lee}
\author[Sewook Oh]{Sewook Oh}

\address{School of Mathematics, Korea Institute for Advanced Study, 85 Hoegiro Dongdaemun-gu, Seoul 02455, Republic of Korea}
\email{juyounglee@kias.re.kr}

\address{Department of Mathematical Sciences and RIM,  Seoul National University, Seoul 08826, Republic of Korea}
\email{shklee@snu.ac.kr}

\address{Center for Mathematical Challenges, Korea Institute for Advanced Study, 85 Hoegiro Dongdaemun-gu, Seoul 02455, Republic of Korea}
\email{sewookoh@kias.re.kr}

\begin{document} 

\begin{abstract} We study  the elliptic maximal functions defined by averages over  ellipses and rotated ellipses which are multi-parametric variants of the circular maximal function. We prove that 
those maximal functions are bounded on $L^p$ for some $p\neq \infty$. For this purpose, we obtain some sharp multi-parameter local smoothing estimates.
 \end{abstract}
 
 \maketitle

\section{Introduction}
Let $\sigma$ denote the normalized measure on the unit circle $\mathbb{S}^1$, and  for $\theta\in \mathbb T:=[0, 2\pi)$ let $R_\theta$ be the rotation in $\mathbb R^2$  such that $R_\theta(1, 0)=(\cos\theta, \sin\theta)$.  For $\theta\in \mathbb T$ and $t,s\in  \mathbb R_+:=(0,\infty)$,  we denote by   $\dsig$  the measure  on 
the rotated ellipse $\mathbb E_{t,s}^\theta:=\{R_\theta(t\cos u, s\sin u):   u  \in \mathbb T\}$ which is given by
\[  \dsig(f) =\int_{\mathbb{S}^1}f\big(R_\theta(ty_1,sy_2)\big)d\sigma(y).\] 
Let $\mathbb I=[1,2]$. 
We consider the maximal operator 
\[\mathfrak M f(x)=\sup_{(\theta,t,s)\in \mathbb T\times \mathbb I^2}  |f\ast \dsig(x)|, \] 
which was called the {\it elliptic maximal function} in \cite{E}. It is a multi-parametric generalization of the circular function whose $L^p$ boundedness on the optimal range of 
$p$ ($p>2$) was proved by Bourgain  \cite{B}  (\cite{S2}). 
Mapping  property of $\mathfrak M$   was studied by Erdo\v gan \cite{E}, who showed that $\mathfrak M$ is bounded from the Sobolev space  $L^{4}_{\alpha}(\mathbb R^2)$ to $L^4(\R^2)$ if $\alpha>1/6$. 
However, the question of whether $\mathfrak M$ admits a nontrivial $L^p$ ($p\neq \infty$)  bound  has remained unknown before this work (see \cite{C, MR, MRZ, H, LL} for previous works on multiparameter maximal functions associated with surfaces). 

 The main object of this article is to prove 
the following. 

\begin{thm}\label{3para thm}  For $p>12$, there is a constant $C$ such that 
   \Be\label{e:max0} 
   \| \M f\|_{L^p(\mathbb R^2)}\le C\|f\|_{L^p(\mathbb R^2)}.
   \Ee
 \end{thm}

However, it was shown  in \cite{E}  that \eqref{e:max0}  fails if $p\le 4$. 
The optimal range of $p$ for which \eqref{e:max} holds remains open. 
As a consequence of the maximal estimate \eqref{e:max0}  one can deduce some measure theoretical results  concerning collections of 
the rotated ellipses (see, for example, \cite{Ma}). 
 In analogue to  the results concerning the circular maximal function \cite{Schlag, SS, L}, 
$L^p$ improving property of $\mathfrak M$ is also of interest. Using the estimates in what follows, one can easily see that 
$\mathfrak M$ is bounded from $L^p$ to $L^q$ for some $p< q$. However, we do not pursue the matter here.  We will address the issue elsewhere.

Let $\mathbb J\subset \mathbb R_+$ be an  interval. The maximal estimate \eqref{e:max0} is a local estimate; however, it is possible to obtain a corresponding global one. 
Using the frequency localized maximal estimate in this paper and a standard argument relying on the Littlewood--Paley decomposition (\cite{B, Schlag}), one can also show that 
the maximal operator 
\[ \bar{\mathfrak M}  f(x)=
\sup_{{(\theta,t,s)\in \mathbb T\times \mathbb R_+^2:\, t/s\in\mathbb J}}|f\ast \dsig(x)|\]
is bounded on $L^p$ for $p>12$ provided that $\mathbb J$ is a compact subset of $\mathbb R_+$. However, as eccentricity of the ellipse  $\mathbb E_{t,s}^\theta$ increases,  
 $\mathbb E_{t,s}^\theta$ gets close to a line.  Using Besicovitch's construction  and taking rotation into account,  it is easy to see  
that  $L^p$ bound fails for any $p\neq \infty$ if $\mathbb J$ is unbounded or the closure  of $\mathbb J$ contains the zero.

We now consider a 2-parameter maximal operator 
\[ \mathcal  Mf(x)=\sup_{(t,s)\in \mathbb R_+^2} |f\ast \sigma_{\!t,s}^{0}(x)|,  \]
which is a circular analogue of the strong maximal function known to be bounded on $L^p$ for  $p>1$. One may call 
$\mathcal M$ the {\it strong circular maximal} operator.  The next theorem shows existence of  nontrivial $L^p$ bound on $\mathcal M$. 
As far as the authors are aware, no such result has been known before.

\begin{thm}\label{2para thm}
 For $p>4$, there is a constant $C$ such that 
   \Be\label{e:max} 
   \| \mathcal M f\|_{L^p(\mathbb R^2)}\le C\|f\|_{L^p(\mathbb R^2)}.
   \Ee
\end{thm}

A modification of the argument in \cite{E} shows that \eqref{e:max} fails if $p\le 3$. We do not know whether \eqref{e:max} holds for $3<p\le 4$. It is also possible to show analogous results in higher dimensions, that is to say, $L^p$ bound on the strong spherical maximal functions \cite{LLO}. 

\subsubsection*{Multi-parameter local smoothing} 
 As is well known in the study of the circular maximal function, its $L^p$ maximal bounds can be obtained by making use of the local smoothing estimate for the averaging operator $f \mapsto f \ast d\sigma_{t,t}^0$ (\cite{MSS, SS}). Following a similar strategy, one may attempt to combine the Sobolev embedding (see, for example, Lemma \ref{sobo3} below) and the (1-parameter) sharp local smoothing estimate for the $2$-d wave operator  due to Guth--Wang--Zhang \cite{GWZ} (also, see \cite{W}) of order $2/p - \epsilon$ for $p \ge 4$, from which one can easily deduce the estimate
\[ 
 \|f\ast d\sigma_{s,t}^0\|_{L_{x,s,t}^p(\mathbb R^2\times \mathbb I^2)}\lesssim \|f\|_{L^p_{\epsilon-2/p}} 
\]
   for  $p\ge 4$ (see \eqref{st4} below). However, this only yields  $L^p_\epsilon$--$L^p$ estimate for $\mathcal M$.   
 In perspective of the Sobolev imbedding, to prove $L^p$ estimate \eqref{e:max}, we need a slightly stronger smoothing estimate, that is to say,  of order bigger than $2/p$. As to be seen later, this is possible by capturing smoothing phenomena associated with two parameter averages.

Our proofs of Theorems \ref{3para thm} and \ref{2para thm} rely, in fact, on multi-parameter local smoothing estimates \eqref{2para-s} and \eqref{3para-s} in Theorem \ref{smoothing} below, which provide smoothing of sharp orders depending on the number of parameters.
 As far as the authors are aware, no such sharp smoothing estimate has appeared in the literature before. Although we do not need to prove the local smoothing estimates up to the sharp order for the maximal estimate, determining the optimal range of $p$ for those sharp smoothing estimates remains interesting problems.

For $\xi\in \mathbb R^2$ and $(t,s)\in \mathbb R_+^2$, let  $\xi_{t,s} =  (t\xi_1, s\xi_2)$
and 
\[\Phi^\theta_\pm(x,t,s,\xi)=x\cdot\xi\pm|(R_{\theta}^\ast \xi)_{t,s}|.\]
Here, $R_{\theta}^\ast$ denotes the transpose of $R_{\theta}$.  Let $B(x,r)$ denote the ball centered at $x$ with radius $r$. 
The asymptotic expansion of the Fourier transform  of ${d\sigma}$ (see \eqref{bessel} below)  naturally leads us to consider the  operators 
\Be
\label{udef} \mathcal{U}^\theta_\pm f(x,t,s)=a(x,t,s)\int e^{i \Phi^\theta_\pm(x,t,s,\xi)}\widehat{f}(\xi) d\xi, 
\Ee
where  $a\in \mathrm C_c^\infty(B(0,2)\times (2^{-1}, 2^2)\times (2^{-1}, 2^2))$.  
The following are our main estimates, which play crucial roles in proving estimates \eqref{e:max0} and \eqref{e:max}.

\begin{thm}\label{smoothing}
If $p\ge 12$ and $\alpha>1/2-3/p$, then   the estimate 
\Be
\label{2para-s}
\| \mathcal{U}^0_\pm f\|_{L^p_{x,t,s}}\le C\|f\|_{L^p_{\alpha}}
\Ee
holds. 
Let $\Delta= \{(t,s)\in  (2^{-1}, 2^2)^2: s=t\}$. Suppose that $\supp a(x,\cdot)\cap \Delta=\emptyset$ for all $x\in B(0,2)$. 
Then, if $p\ge 20$  and $\alpha> 1/2-4/p$,  we have  the estimate 
\Be 
\label{3para-s}
\|  \mathcal{U}^\theta_\pm f\|_{L^p_{x,t,s,\theta}}\le C\|f\|_{L^p_{\alpha}}. 
\Ee 
\end{thm}

Compared with the local smoothing estimate for the 2-d wave operator $f\mapsto \mathcal{U}^0_\pm f(\cdot,t,t)$, the estimates  \eqref{2para-s} and \eqref{3para-s} have 
additional  smoothing of order up to $1/p$ and $2/p$, respectively, which results from averages in $s,t$;  and  $s,t$, and $\theta$. 
The smoothing orders in  \eqref{2para-s} and \eqref{3para-s} are sharp (see Section \ref{remarks}) in that \eqref{2para-s} and \eqref{3para-s} fail if $\alpha<1/2-3/p$ and $\alpha<1/2-4/p$, respectively.  The ranges of $p$ for  which \eqref{2para-s} and \eqref{3para-s} hold are actually determined  by our use of the decoupling inequalities  (see  Theorem  \ref{BDGcone}  and  \ref{var decoupling}). However, there is no reason to believe that  those  ranges of $p$ are optimal.    It is clear that  the condition $\supp a(x,\cdot)\cap \Delta=\emptyset$  is necessary for \eqref{3para-s} to hold with all $\alpha>1/2-4/p$. 
Indeed, when $t$ get close to $s$,  the ellipse $\mathbb E_{t,s}^\theta$ becomes close to the circle $\mathbb E_{s,s}^\theta$, which is invariant under rotation,  so that  average in $\theta$ does not 
yield any further regularity gain.

An immediate consequence of the estimate \eqref{2para-s} is that the two-parameter averaging operator $f\mapsto a (f\ast \dsi)$ is bounded from 
$L^p_{x,t,s}$ to $L^p_{\alpha}$  for $\alpha<3/p$.  
From those estimates with $p\neq \infty$, following the argument in \cite{HKLO},  one can also obtain results regarding dimensions of unions of ellipses.

In \cite{LL},   the first and the second authors studied 2-parameter local smoothing estimates for  averages over tori, $t\mathbb S^1\times s\mathbb S^1$ while $t,s$ are contained in a compact set.  As a consequence, they  completely characterized $L^p$--$L^q$ estimates for the associated 2-parameter maximal function.  Their results  were  obtained  by making use of  the (one-parameter) sharp local smoothing estimate for the 2-d wave operator  (\cite{GWZ}).  The main geometric observation  was that  a torus can essentially be  regarded as a product of two circles. This idea becomes more transparent  if one considers the 2-parameter propagator
\[   \int_{\R^3} e^{i(x\cdot \eta+ t\sqrt{\eta_1^2+\eta_2^2}+s|\eta|)}\widehat{g}(\eta)d\eta, \]
which constitutes the major part of the averaging operators. Roughly speaking, by comparing size of $\sqrt{\eta_1^2+\eta_2^2}$  and $|\eta_3|$ in the frequency side,  one can break the phase function into two independent ones with respect to two parameters, $t,s$, so that  the sharp 1-parameter local smoothing estimate (see Lemma \ref{lem:locals} below)  can be applied by fixing some variables. 

  In contrast, the multiparameter averages considered in this article do not exhibit  such nice properties.  Neither of  the phase functions of  the operators $\mathcal{U}_{\pm}^{0}$, $\mathcal{U}_{\pm}^{\theta}$ is linear in $(x,t,s)$, $(x,t,s,\theta)$,  respectively.   Moreover,  it is highly unlikely  that one can obtain  the sharp local smoothing estimates for $\mathcal{U}_{\pm}^{\theta}$ and $\mathcal{U}_{\pm}^{0}$ by making use of those for 1-parameter operators.   This compels us  to consider the muliti-parameter local smoothing from a more genuine perspective.


\subsubsection*{Key observation}  Now we briefly explain the key  geometric properties which allow us  to obtain the sharp local smoothing estimates.  The main ingredients for the proof of the  estimates \eqref{2para-s} and \eqref{3para-s} are decoupling inequalities for the operators $\mathcal{U}^0_\pm$ and $\mathcal{U}^\theta_\pm$ 
(see, for example, Lemma \ref{2paracase} and \ref{0->4} below).  
Those inequalities are built  on our striking observation  that   the immersions 
\begin{align}
& 
\label{th0}
 \xi\mapsto \nabla_{x,t,s}  \Phi^0_\pm(x,t,s,\xi),
 \\
 &
\label{th} 
\xi\mapsto \nabla_{x,t,s,\theta}  \Phi^\theta_\pm(x,t,s,\xi)
 \end{align}
(fixing  $(x,t,s)$ and $(x,t,s,\theta)$ with $s\neq t$, respectively) 
give rise to submanifolds, which are conical extensions of a  finite type curve in $\mathbb R^3$ and a nondegenerate  curve in $\mathbb R^4$,   respectively.  
By this observation, we are naturally  led  to regard the operators $\mathcal{U}^\theta_\pm, \mathcal{U}^0_\pm$ as  variable coefficient generalizations of the operators given by the conic surfaces. 

 Meanwhile, the decoupling inequalities for the  extension operators given by those conic surfaces,  which are constant coefficient counterparts of the abovementioned operators,  are already known (see Theorem \ref{BDGcone} below and \cite{BGHS2}). Those inequalities were, in fact, deduced from the decoupling inequality for the nondegenerate curve  due to Bourgain--Demeter--Guth \cite{BDG}. To obtain such  inequalities for  $\mathcal{U}^0_\pm$ and $\mathcal{U}^\theta_\pm$, we combine the known inequalities for the extension (adjoint restriction)  operators and the argument in   \cite{BHS} to get desired  decoupling inequalities in a variable coefficient  setting (see Theorem \ref{var decoupling}).

Decoupling inequalities of different forms have been extensively  used in the recent studies on maximal and smoothing estimates for 
averaging operators. We refer the reader to \cite{PS, KLO2, KLO,  BGHS} and references therein for related works. 

\smallskip

\noindent{\bf Organization.}
In Section \ref{multipara section}, we prove  Theorem \ref{3para thm} and \ref{2para thm} by mainly making use of  the multi-parameter smoothing estimates in Theorem \ref{smoothing} (see also  Proposition \ref{main u2} and \ref{main u3} below). In Section \ref{decoup sec}, we discuss a variable coefficient generalization of the decoupling inequality for the aforementioned conic surfaces.  In Section \ref{nondegeneracy section}, we prove our main estimates in Proposition \ref{main u2} and \ref{main u3} making use of the variable coefficient decoupling inequality, of  which proof we provide in Appendix.

\smallskip

\noindent{\bf Notation.} Let 
$\mathbb B^{k}(x,r)$ denote the ball in $\mathbb R^{k}$ which is centered at $x$  and of radius $r$. However, if there is no ambiguity  concerning the dimension of the ambient space,  we simply denote by $B(x,r)$ the ball centered at $x$ with radius $r$ as above. 
In addition to $\,\widehat{}\,\,$ we also use $\mathcal F$ to denote the Fourier transform. 

\smallskip

\noindent{\bf Added in proof.}  Shortly after this paper was uploaded to arXiv, a related work by Chen, Guo, and Yang \cite{CGY} appeared, where they generalized the local multiparameter smoothing estimates under a cinematic curvature condition. The $L^p$   maximal bounds in Theorems \ref{2para thm} and \ref{3para thm} are now established on the optimal ranges by utilizing the estimates in \cite{PYZ} and \cite{Zahl}, along with the local smoothing estimates in Propositions \ref{main u2} and \ref{main u3} (see also \cite{CGY}).

\section{Proof of maximal bounds}\label{multipara section}

In this section we prove the maximal estimates while assuming the local smoothing estimates. We begin by recalling  an elementary lemma, which is easy to show using the fundamental theorem of calculus and H\"older's inequality (for example, see \cite[Lemma 3.1]{LL}).

\begin{lem}
\label{sobo3}
Let $1\le p\leq \infty$, and  $J_1$, $J_2$, and $J_3$ be closed intervals of length $\sim 1$. 
Let $\mathfrak R=J_1\times J_2\times J_3$ and $G\in \mathrm C^1( \mathfrak R)$.  Then, there is a constant $C>0$ such that 
\begin{align*}
 & \sup_{(t,s,\theta)\in \mathfrak R}\vert G(t,s,\theta) \vert \le C (\lambda_1\lambda_2\lambda_3)^{\frac{1}{p}}\sum_{\beta\in\{0,1\}^3}(\lambda_1^{-1},\lambda_2^{-1},\lambda_3^{-1})^{\beta}
        \Vert \partial_{t,s,\theta}^\beta G\Vert_{L^p(\mathfrak R)} 
 \end{align*}
holds for  $\lambda_1,\lambda_2,\lambda_3\ge 1$. Here  $\beta=(\beta_1, \beta_2, \beta_3)$ denotes a triple index. 
\end{lem}

By the Fourier inversion we write
\Be 
\label{ft}
 f\ast \dsig(x)=(2\pi)^{-2}\int e^{ix\cdot\xi}\,\widehat{f}(\xi)\,\widehat{d\sigma}\big((R_{\theta}\xi)_{t,s}^\ast\xi\big)\,d\xi. 
 \Ee
We recall the well known asymptotic expansion  of $\widehat\sigma$ (for example, see \cite{S}):
\begin{equation}\label{bessel}
  \textstyle   \widehat{\sigma}(\xi)=\sum_{\pm, \,0\leq k\leq N}\,C_j^{\pm}|\xi|^{-\frac{1}{2}-k}e^{\pm i|\xi|}+E_N(|\xi|),\quad |\xi|\geq 1,
\end{equation}
where $E_N$ is a smooth function satisfying
$
    |({d}/{dr})^l E_N(r)|\leq Cr^{-l-(N+1)/{4}},$ $0\leq l\leq (N+1)/{4}
$
for $r\geq 1$ and a constant $C>0$.  Fixing a sufficiently large $N$, 
we may disregard $E_N$.  Thus, it suffices to consider the contribution from the main part $k=0$ since the remaining parts can be handled similarly but more easily. 

The following is an easy consequence of the sharp local smoothing estimate for the wave operator 
 \[ \mathcal{W}_{\pm}  f (x,t)\coloneqq \int e^{i(x\cdot\xi\pm t|\xi|)} \widehat{f}(\xi)d\xi.\]

\begin{lem}[{\cite[Lemma 2.2]{LL}}]
\label{lem:locals} 
Let $p\ge 4$, $j\ge 0$, $ b\gtrsim 2^{-j/2}$ and $\mathbb{A}_j:=\{ \xi\in\R^2 : 2^{j-1}\leq |\xi|\leq 2^{j+1}\}$. Then, for any $\epsilon>0$  there is a constant $C_\epsilon$ such that 
\Be \label{locals}
\big\|
\mathcal W_\pm f
\big\|_{L^p (\R^2\times [0,2])}
\le C_\epsilon
 b^{1-\frac4p} 2^{(\frac12-\frac2p +\epsilon)j}  \Vert f\Vert_{L^p}
\Ee
whenever $\supp \widehat f\subset \mathbb{A}_j$ and  is included in an angular sector of angle  $b$.
\end{lem}

Indeed, this follows from interpolation between  the estimate \eqref{locals} for $p=\infty$ and for $p=4$. The latter is basically due to Guth, Wang, and Zhang \cite{GWZ}. 
Using this estimate, we can get the following estimate:
\Be 
\label{st4}
\Vert \mathcal U^\theta_\pm f\Vert_{L^4_{x,t,s} } \lesssim 2^{\epsilon j}\Vert f\Vert_{L^4} 
\Ee
for any $\epsilon>0$ provided that  $\supp \widehat f\subset \mathbb{A}_j$. 
We will employ \eqref{st4} and its variant  to obtain some (non-sharp) smoothing estimates on an extended range by interpolation (see Proposition \ref{main u2} and \ref{main u3} below). 
 Indeed, note that 
    \Be 
    \label{utheta}   \tilde{\mathcal U}^{\theta,s}_\pm  f(x,t):=   \mathcal U^\theta_\pm  f (x,t, ts)=a(x,t,ts)\int e^{i(x\cdot\xi\pm t|(R_{\theta}^\ast \xi)_{1,s}|)}\widehat{f} (\xi)d\xi. 
    \Ee
By a change of variables  and the $L^4$ local smoothing estimate for $\mathcal W_\pm$, we have  
$ \Vert \tilde{\mathcal U}^{\theta,s}_\pm f\Vert_{L^4_{x,t} } \le C 2^{\epsilon j}\Vert f\Vert_{L^4} $ for any $(\theta, s)$ whenever $\supp \widehat f\subset \mathbb{A}_j$.
Taking integration in $s$, we get \eqref{st4}.

\subsection{2-parameter maximal operator $\mathcal M$: Proof of Theorem \ref{2para thm}
}  To show Theorem \ref{2para thm},  we make use of the following, which we prove in Section \ref{nondegeneracy section}.

\begin{prop}\label{main u2}
    Let $4\leq p\leq \infty$.    For any $\epsilon>0$, we have 
    \begin{equation}
    \label{2-local}
        \Vert \mathcal{U}^0_\pm f\Vert_{L^p_{x,t,s}}\lesssim \begin{cases}
            2^{(\frac{3}{8}-\frac{3}{2p}+\epsilon)j}\Vert f\Vert_{L^p}, &\,\, \  4\leq p <12,\\
            2^{(\frac{1}{2}-\frac{3}{p}+\epsilon)j}\Vert f\Vert_{L^p},  &  \ 12\le p\leq \infty
        \end{cases}
    \end{equation}
    whenever  $\supp \widehat f\subset \mathbb{A}_j$.
\end{prop}

Note that the estimate \eqref{2-local} is equivalent to \eqref{2para-s} in Theorem \ref{smoothing} when $12\le p\le \infty$. As mentioned earlier, this range of $p$ is subject to that of  Theorem \ref{var decoupling} with $d=3$.  For $4\leq p< 12$, \eqref{2-local} follows from interpolation between  \eqref{2-local} with $p=12$ and \eqref{st4}. 
 Instead, one may try to interpolate the estimate  with the trivial $L^2$ estimate  but this only yields weaker estimates.

  To prove \eqref{e:max},  we  consider a local maximal operator 
\[ \mathcal  M_{loc} f(x)=\sup_{(t,s)\in (0,2]^2} |f\ast \sigma_{\!t,s}^{0}(x)|.  \]
By scaling it is sufficient for the proof of Theorem \ref{2para thm} to show  
\Be 
\label{e:maxc}
 \| \mathcal M_{loc} f\|_{L^p(\mathbb R^2)}\le C\|f\|_{L^p(\mathbb R^2)}, \quad p>4. 
\Ee

For the purpose, we extensively use the Littlewood-Paley projection to decompose $f$. 
Let  $\phi$ denote a smooth function such that $\supp \phi\subset (1-2^{-10},2+2^{-10})$ and $\sum_{j=-\infty}^{\infty}\phi(\tau/2^j)=1$ for $\tau>0$. 
For simplicity, we set  $\phi_j(\tau)=\phi(\tau/2^j)$ and 
\[ \textstyle \phi_{<k}(\tau)=\sum_{j<k}\phi_j(\tau). \]
We define $\phi_{\geq k}$ similarly.  For $(k,n)\in\Z^2$, let  
\begin{align*}
 \widehat{f_k^n}(\xi_1, \xi_2)&=  \phi_k(|\xi_1|)\phi_n(|\xi_2|) \widehat f(\xi_1, \xi_2), 
 \\
  \mathcal F(f_{<k}^{<n})(\xi_1, \xi_2)&=  \phi_{<k}(|\xi_1|)\phi_{<n}(|\xi_2|) \widehat f(\xi_1, \xi_2).
  \end{align*}
 Similarly, $f_{<k}^{\ge n}$, $f_{\ge k}^{<n}$,  and  $f_{\ge k}^{\ge n}$ are defined. In particular, for any $n,k\in \mathbb Z$, we have 
\Be\label{fff}
f=f_{<k}^{<n} + f_{\geq k}^{<n}+f_{<k}^{\geq n}+ f_{\geq k}^{\geq n}.\
\Ee

\begin{proof}[Proof of \eqref{e:maxc}] 
Denoting $Q_{k}^n=[2^{-k},2^{-k+1}]\times [2^{-n},2^{-n+1}]$, 
we set 
\begin{align*}
     & \mathcal{M}_1f=\sup_{k,n\geq 0}\sup_{(t,s)\in Q_{k}^n}|f_{<k}^{<n} \ast \dsi|, 
     \\& \mathcal{M}_2f=\sup_{k,n\geq 0}\sup_{(t,s)\in Q_{k}^n}|f_{\geq k}^{<n}\ast \dsi|,
     \\
     & \mathcal{M}_3f=\sup_{k,n\geq 0}\sup_{(t,s)\in Q_{k}^n}|f_{<k}^{\geq n}\ast \dsi|, 
     \\
     & \mathcal{M}_4f=\sup_{k,n\geq 0}\sup_{(t,s)\in Q_{k}^n}|f_{\geq k}^{\geq n}\ast \dsi|.
\end{align*}
Since $ \mathcal M_{loc} f(x)=\sup_{k,n\geq 0}\sup_{(t,s)\in Q_{k}^n} |f\ast \sigma_{\!t,s}^{0}(x)|$, from \eqref{fff}  it follows that 
\[\textstyle \mathcal M_{loc} f(x)\leq \sum_{j=1}^4 \mathcal{M}_j f(x)
. \]
The maximal operators
$\mathcal{M}_1$, $\mathcal{M}_2$, and $\mathcal{M}_3$ can be handled easily. We note that $f_{<k}^{<n} =f\ast  K$ with a kernel $K$ satisfying  
\[  |K(x)|\lesssim   \mathrm K_k^n(x):={2^{k+n}}{(1+2^k|x_1|)^{-N}(1+2^n|x_2|)^{-N}} \]
for any large  $N$. Hence, it follows that $ |f_{<k}^{<n}\ast \sigma_{t,s}^0(x)|  \lesssim  \mathrm K_k^n\ast |f|(x)$ if $(t,s)\in Q_k^n$.   This gives 
$   \mathcal M_1f(x)     \lesssim M_sf(x)$ where $M_s$ denotes the strong maximal operator. Therefore, we get $\| \mathcal M_1f\|_p\lesssim \|f\|_p$ for $1<p\le \infty$. 

 We denote  by $H$ the one dimensional Hardy-Littlewood maximal operator. 
 For the maximal operator $\mathcal {M}_2$, note that  $\mathcal F(f_{\geq k}^{<n})=\widehat f(\xi) \phi_{<n}(|\xi_2|)-  \mathcal F(f_{<k}^{<n})(\xi)$. 
Thus, as before,  we observe that 
\[ |f_{\geq k}^{<n}\ast \dsi(x)|\lesssim \iint \frac{2^n| f (x_1-ty_1,z_2)|}{(1+2^n|x_2-z_2|)^N} d\sigma(y)dz_2  + \mathcal{M}_1f(x) \]
for $s\sim 2^{-n}$. This yields 
\[ \mathcal M_2f(x)\lesssim   H(M_cf (x_1, \cdot))(x_2)+ \mathcal{M}_1f(x),\]
where $M_c h(x_1)=\sup_{0<t< 2} \int  h(x_1-ty_1) d\sigma(y)$.  
 Using Bourgain's circular maximal theorem it is easy to  see that $M_c$ is bounded on $L^p(\mathbb R)$ for $p>2$ (see \cite[Lemma 3.2]{LL}). Consequently, $L^p$ boundedness of $M_s$ and $H$ yields  
$ \Vert \mathcal{M}_2 f\Vert_{p}\lesssim \Vert f\Vert_{p} $ for $p>2$.  A symmetric argument also shows that  $ \Vert \mathcal{M}_3 f\Vert_{p}\lesssim \Vert f\Vert_{p} $
for $p>2$. 

Finally, we consider $\mathcal M_4f$, which constitutes the main part.  We note that $\mathcal M_4f\le \sup_{k,n\geq 0}\sum_{j,l\geq 0}\sup_{(t,s)\in Q_{k}^n}|f_{j+k}^{l+n}\ast d\sigma^0_{s,t}|$. 
The imbedding $\ell^p\hookrightarrow \ell^\infty$, followed by Minkowski's inequality, gives
\begin{align*}
   \Vert \mathcal M_4f \Vert_{p} 
    &\leq \sum_{j,l\geq 0}\Big(\sum_{k,n\geq 0} \Big\| \sup_{(t,s)\in Q_k^n}\big|f_{j+k}^{l+n}\ast \dsi\big|  \Big \Vert_{p}^p\,\Big)^{1/{p}}.
\end{align*}
We now claim that
\Be
\label{jj}
\Big\Vert \sup_{(t,s)\in Q_k^n}\big|f_{j+k}^{l+n}\ast \dsi\big| \Big \Vert_{p}\!\lesssim 2^{-\delta\max(j, l)}\|f_{j+k}^{l+n}\|_p, \quad j, l\ge 0
\Ee
for some $\delta>0$ if $p>4$. Once we have this, it is easy to show that $\mathcal M_4$ is bounded on $L^p$ for $p>4$. Indeed, note that  $(\sum_{k,n} \|f_{j+k}^{l+n}\|_{L^p}^p)^{{1}/{p}}\le C \|f\|_p$ for $2\le p\le \infty$, which follows by interpolation between the estimates for $p=2$  and $p=\infty$ (see, for example, \cite[
Lemma 6.1]{TVV}).  Combining \eqref{jj}  and  this inequality  gives
\[ \textstyle  \Vert \mathcal M_4f \Vert_{L^p} 
\lesssim    \sum_{j,l\geq 0}  2^{-\delta\max(j, l)}\|f\|_p\lesssim \|f\|_p.  \]

 It remains to show \eqref{jj}.  By rescaling we note that the two operators  $f\mapsto \sup_{(t,s)\in Q_0^0}|f_j^{l}\ast \dsi|$ and  $f\mapsto \sup_{(t,s)\in Q_k^n}|f_{j+k}^{l+n}\ast \dsi|$ have the same bounds on $L^p$. Thus, it suffices to show \eqref{jj} for $k=n=0$.  To this end, by the finite speed of propagation and translation invariance, it is enough to prove that
\[
\Big\Vert \sup_{(t,s)\in Q_0^0}\big|f_j^{l}\ast \dsi\big| \Big \Vert_{L^p(B(0,1))}\lesssim 2^{-\delta\max(j, l)}\|f_j^{l}\|_p, \quad j, l\ge 0.
\]
Note that the Fourier support of $f_j^{l}$ is included in $\mathbb{A}_{\max(j,l)}$.
We recall  \eqref{ft} and  \eqref{bessel}. So,  it is enough to consider, instead of $f\to f_j^{l}\ast \dsi$,  the operators 
\[
  \mathcal A_\pm f(x,t,s)=  \int |\xi_{t,s}|^{-\frac{1}{2}}e^{i(x\cdot\xi\pm |\xi_{t,s}|)} \widehat{f_j^{l}}(\xi)\, d\xi.
\]
Contributions from other terms in \eqref{bessel} can be handled similarly but they are less significant. 
Therefore,   the matter is reduced to obtaining the estimate
\Be
\label{jj2}
\Big\Vert \sup_{(t,s)\in Q_0^0}\big|\mathcal A_\pm f\big| \Big \Vert_{L^p(B(0,1))}\lesssim 2^{-\delta\max(j, l)}\|f_j^{l}\|_p, \quad j, l\ge 0. 
\Ee
Since $\partial_t|\xi_{t,s}|={t\xi_1^2}/{|\xi_{t,s}|}$ and $\partial_s|\xi_{t,s}|={s\xi_2^2}/{|\xi_{t,s}|}$,  applying  Lemma \ref{sobo3} (with $\lambda_1=\lambda_2=2^{\max(j, l)}$ and 
$\lambda_3=1$) to  $\mathcal A_\pm f$  and using Mikhlin's multiplier theorem, 
we get  
\[  \Big\Vert \sup_{(t,s)\in Q_0^0}\big|\mathcal A_\pm f| \Big \Vert_{L^p(B(0,1))}\lesssim  2^{(\frac2p-\frac12)\max(j, l)}   
 \big\Vert  \mathcal U^0_\pm f_j^{l} \big \Vert_{L^p_{x,t,s}}. \]
 Since the Fourier support of $f_j^{l}$ is included in $\mathbb{A}_{\max(j,l)}$, by Proposition \ref{main u2}  it follows that \eqref{jj2} 
holds for some $\delta>0$ as long as $p>4$.     \end{proof}

\subsection{3-parameter maximal operator $\mathfrak M$: Proof of Theorem \ref{3para thm}}
The proof  basically relies on   the estimate \eqref{3para-s}.  However, to control the  averages when  $s,t$ are close to each other, we need an additional decomposition: 
\Be 
\label{udef1} \textstyle \mathcal U^\theta_\pm f(x,t,s) =  \sum_k \mathcal U^\theta_{\pm, k} f(x,t,s):=  \sum_k \psi_k(t,s) \mathcal U^\theta_\pm f(x,t,s),
\Ee
where $\psi_k(t,s)= \phi(2^k|s-t|).$  Note that $\mathcal U^\theta_{\pm, k}=0$ if $k\le -3$.

\begin{prop}\label{main u3}
    Let $4\leq p\leq \infty$ and $0\le k\le j$. For any $\epsilon>0$, we have
        \begin{equation}\label{main u3 estimate}
        \Vert  \mathcal U^\theta_{\pm, k} f\Vert_{L^p_{x,t,s,\theta}}\lesssim \begin{cases}
            2^{(\frac{3}{8}-\frac{3}{2p}+\epsilon)j}2^{\frac{k}{p}}\Vert f\Vert_{L^p}, & \,\,4\leq p< 20,\\
            2^{(\frac{1}{2}-\frac{4}{p}+\epsilon)j}2^{\frac{k}{p}}\Vert f\Vert_{L^p}, & 20\le p\leq \infty 
        \end{cases}
    \end{equation}
 whenever   $\supp \widehat f\in \mathbb{A}_j$. 
\end{prop}

We also note that  \eqref{main u3 estimate} is equivalent to the estimate  \eqref{3para-s} in Theorem \ref{smoothing} for $20\le p\leq \infty$. As before, the range of $p$ corresponds to that of  Theorem \ref{var decoupling} with $d=4$.  For $4\leq p< 20$,  \eqref{main u3 estimate} follows from \eqref{st4} and  \eqref{main u3 estimate} with $p=20$ via interpolation.   

Once we have Proposition \ref{main u3},  
the proof of Theorem \ref{3para thm} proceeds in a similar  manner as  that of Theorem \ref{2para thm}.      
    Let $f_j$, $f_{<j}$ be, respectively,  defined by  
\[ \textstyle \widehat{f_j}(\xi)=\phi_j(|\xi|)\widehat{f}(\xi), \quad \widehat{f_{<j}}(\xi)=\sum_{l<j} \phi_{l}(|\xi|)\widehat{f}(\xi). \] 
Note  that $\sup_{(\theta,t,s)\in   \mathbb T\times \mathbb I^2}  |f_{<1} \ast \dsig(x)|\lesssim K_N\ast |f|(x) $ for any $N$ where  $K_N(x):= (1+|x|)^{-N}$.  
Thus,  it suffices to consider $f\mapsto \tilde {\mathfrak M} f:= \sum_{j\ge 1}\sup_{(\theta,t,s)\in   \mathbb T\times \mathbb I^2}   | f_j\ast \dsig|$. 
We  make  decomposition in  $t,s$ using $\psi_k$ to get 
\[
\textstyle \tilde {\mathfrak M} f \le   \mathfrak M' f+\mathfrak M'' f:= \sum_{j\ge 1} \sum_{ k\le j}   \mathfrak M_{k} f_j+ \sum_{j\ge 1} \sup_{k> j} \mathfrak M_{k} f_j,
\]
where 
\begin{align*}
\mathfrak M_{k} f(x)=   \sup_{(\theta,t,s)\in   \mathbb T\times \mathbb I^2}  |\psi_{k}(t,s)\, f\ast \dsig(x)|. 
\end{align*}
The operator    
$\mathfrak M''$ can be handled by using the  bound on the circular maximal function. Indeed, observe that
\[ 2^{2j} K_N(2^j\cdot)\ast \dsig(x)\lesssim 2^j (1+ 2^j||(R^\ast_\theta x)_{1,t/s}|-t|)^{-N+2}\] 
for $t,s\in \mathbb I$. This gives $2^{2j}|\psi_{k}(t,s)| K_N(2^j\cdot)\ast \dsig(x)
\lesssim 2^j (1+ 2^j||x|-t|)^{-N+2}$ for $k\ge j$ because $|t-s|\lesssim 2^{-j}$.  Note $|f_j\ast \dsig|\lesssim |f_j|\ast 2^{2j} K_N(2^j\cdot)\ast \dsig$. 
So, combining these inequalities and  taking  $N$ sufficiently large, we see that 
\[\textstyle \sup_{k> j} \mathfrak M_{k} f_j\lesssim  \bar {\mathcal M}f_j(x)+2^{-10 j} K_{10}\ast |f_j|(x),\] 
where $\bar {\mathcal M} g(x)=\sup_{t\in (2^{-1}, 2^2)} |g\ast \sigma_{t,t}^0(x)|.$  It is well known that $\|\bar {\mathcal M}f_j\|_p\le 2^{-cj} \|f\|_p$ for some $c>0$ if $p>2$
(see \cite{MSS, L}).
Therefore,  $\mathfrak M''$ is bounded on $L^p$ for $p>2$.

To show $L^p$ bound on $\mathfrak M'$, as before, we only need to obtain a local estimate 
$\|  \mathfrak M' f\|_{L^p(B(0,1))}\lesssim \|f\|_p$ 
for $p>12$. This is immediate once we  have
\[
\Vert  \mathfrak M_{k} f_j\Vert_{L^p(B(0,1))}\lesssim  2^{-\epsilon_0 j}
\Vert f\Vert_{L^p}, \quad 1\le k\le j
  \]
  for any $p>12$ and some $\epsilon_0>0$. 
By  \eqref{ft} and \eqref{bessel}, the estimate follows if we show
 \Be 
    \label{hoho0} \Big\Vert  \sup_{(\theta, t,s)\in \mathbb T\times \mathbb I^2} |\mathcal U^\theta_{\pm ,k} f_j|\Big\|_{p}\lesssim 
    2^{(\frac{3}{8}+\frac{3}{2p}+\epsilon)j}\Vert f\Vert_{p},  \quad 4\leq p\leq 20 
 \Ee
for any $\epsilon>0$. 

\begin{proof}[Proof of \eqref{hoho0}]
    We  use Lemma \ref{sobo3}.  To do so,  we observe  that 
\begin{equation}
\label{comp}    \begin{aligned}
        &\partial_t|(R_{\theta}^\ast\xi)_{t,s}|=m_1:={t(R_{\theta}^\ast\xi)_1^2}{|(R_{\theta}^\ast\xi)_{t,s}|}^{-1},\\
        &\partial_s|(R_{\theta}^\ast\xi)_{t,s}|=m_2:={s(R_{\theta}^\ast\xi)_2^2}{|(R_{\theta}^\ast\xi)_{t,s}|}^{-1},\\
        &\partial_\theta|(R_{\theta}^\ast\xi)_{t,s}|=m_3:=2(t^2-s^2)(R_{\theta}^\ast\xi)_1(R_{\theta}^\ast\xi)_2{|(R_{\theta}^\ast\xi)_{t,s}|}^{-1},
    \end{aligned}
    \end{equation}
    where $R_{\theta}^\ast\xi=(  (R_{\theta}^\ast\xi)_1, (R_{\theta}^\ast\xi)_2 )$. 
   It is clear that $|\partial_\xi^\alpha m_l |\lesssim  |\xi|^{1-|\alpha|}$, $l=1,2$, and $|\partial_\xi^\alpha m_3|\lesssim 2^{-k} |\xi|^{1-|\alpha|}$. 
    Note that $|\nabla_{t,s} \psi_k(t,s)|\lesssim 2^k\le 2^j$. Recalling \eqref{udef} and \eqref{udef1},  we apply Lemma  \ref{sobo3}  to $\sup_{(\theta, t,s)\in \mathbb T\times \mathbb I^2} |\mathcal U^\theta_{\pm ,k} f_j|$    with  $\lambda_1=\lambda_2=2^j$ and $\lambda_3=2^{j-k}$.  
   Thus,    by Mikhlin's multiplier theorem  we have 
  \[ \Big\Vert  \sup_{(\theta, t,s)\in \mathbb T\times \mathbb I^2} |\mathcal U^\theta_{\pm ,k} f_j|\Big\|_{p}\lesssim  2^{(3j-k)/p} \big \Vert  \mathcal U^\theta_{\pm ,k} f_j\big\|_{L^p_{x,t,s,\theta}}. \] 
 By Proposition \ref{main u3} the  estimate \eqref{hoho0} follows. 
\end{proof}

\subsection{Remarks on the estimates \eqref{2para-s} and \eqref{3para-s}}
\label{remarks}
Before closing this section, we make remarks regarding the  estimates \eqref{2para-s} and \eqref{3para-s}. 

$(i)$ Once one has the estimates \eqref{2-local} and \eqref{main u3 estimate}, the proofs of the estimates \eqref{2para-s} and \eqref{3para-s} are rather straightforward. So, we omit them. 

 $(ii)$ As mentioned before, the smoothing orders in  the estimates \eqref{2para-s} and \eqref{3para-s}  are sharp, that is to say,   \eqref{2para-s}, \eqref{3para-s} fail if $\alpha<1/2-3/p$, $\alpha<1/2-4/p$, respectively. To see this, 
we only consider the operator $ \mathcal U^\theta_+ $ since $ \mathcal U^\theta_- $ can be handled similarly. 
Let $g$ be a function given by $ \widehat g(\xi)= \phi(2^{-j}|\xi_{1,2}|) e^{-i|\xi_{1,2}|}.$
It is easy to see that $\|g\|_{L^p_\alpha}\lesssim 2^{(\alpha+3/2-1/p)j}$. Note that 
\[  \mathcal U^\theta_+ g(x,t,s)=2^{2j} \int e^{i2^j( x\cdot \xi+ |(R_\theta^\ast\xi)_{t,s}| -|\xi_{1,2}|)} \phi(|\xi|) d\xi. \] 
Thus, we have $ |\mathcal U_+^0 g(x,t,s)|\gtrsim 2^{2j}$ if $|x|, |t-1|,$ and $|s-2|\le 2^{-j}/100$. Consequently, the estimate \eqref{2para-s} gives   
$2^{(2-4/p)j}\lesssim  2^{(\alpha+3/2-1/p)j} $. Letting $j\to \infty$ shows that \eqref{2para-s} holds only if  $\alpha\ge 1/2-3/p$. Concerning the estimate \eqref{3para-s},  note that 
$|\mathcal U_+^\theta g(x,t,s)|\gtrsim 2^{2j}$ if $|x|, |\theta|,  |t-1|,$ and $|s-2| \le 2^{-j}/100$. So, the estimate \eqref{3para-s} implies
$2^{(2-5/p)j}\lesssim  2^{(\alpha+3/2-1/p)j} $. Therefore,  \eqref{3para-s} holds only if  $\alpha\ge 1/2-4/p$. 

$(iii)$ Besides those lower bounds on the smoothing order $\alpha$, we also have
\Be\label{ad-nec}   \alpha\geq {1}/{4}-{1}/{p}, \Ee
which is valid for both \eqref{2para-s} and \eqref{3para-s}.  However, we are currently far from being able to prove the estimates of smoothing orders up to these bounds, except in some trivial cases. This problem seems to be highly challenging.
 
  To show \eqref{ad-nec}, let us consider  $\widehat{h}(\xi)=e^{-i|\xi|}\phi(2^{-j}\xi_1)\tilde{\phi}(2^{-3j/4}\xi_2)$ where $\tilde{\phi}\in C_c^{\infty}([-2,2])$ such that $\tilde{\phi}(\tau)=1$ if $-1\leq \tau\leq 1$. It is straightforward to verify  that $\| h\Vert_{L^p_{\alpha}}\lesssim 2^{(\alpha+3/2-5/4p)j}$. Similarly to the above argument, we obtain that $|\mathcal{U}_+^0 h (x,t,s)|\gtrsim 2^{7j/4}$ if $-1<x_1<0$, $|x_1+t-1|\lesssim 2^{-j}$, $|{s^2}/{t}-1|\lesssim 2^{-j/2}$, and $|x_2|\lesssim 2^{-3j/4}$. Thus, the estimate \eqref{2para-s} gives $2^{({7}/{4}-{9}/{4p})j}\lesssim 2^{(\alpha+{3}/{2}-{5}/{4p})j}.$   
Letting $j\to\infty$, we get  \eqref{ad-nec} for \eqref{2para-s} to hold.  Using the same function $h$ and taking into account  rotational symmetry,  one can easily show that  \eqref{ad-nec} is also necessary for \eqref{3para-s} to hold. 

\medskip

The rest of the paper is devoted to the proof of Proposition \ref{main u2} and \ref{main u3}.

\section{Variable coefficient  decoupling inequalities}\label{decoup sec}

In this section, we discuss the decoupling inequalities which we need to prove Proposition \ref{main u2} and \ref{main u3}. 

\begin{defn}
Let $I$ be an interval and $\gamma:  I \rightarrow \R^d$ be a smooth curve. We say  $\gamma$ is nondegenerate  if
$\det(\gamma'(u),\cdots,\gamma^{(d)}(u))\neq0$ for all $u\in I.$
\end{defn}

For a curve $\gamma$ defined  on  $\mathbb I_0:=[-1,1]$,  set $\mathfrak C(\gamma)=\big\{ r(1,\gamma(u)):   u\in \mathbb I_0,  \ r\in \mathbb I\big\}$, which we call  the {\it conical extension of $\gamma$}. 
Consider  an adjoint restriction operator
\begin{equation}\label{def of E}
E^\gamma g(z):=\iint_{\I_0\times \mathbb I} e^{iz\cdot r(1,\gamma(u))} g(u,r)dudr,  \quad z\in \mathbb R^{d+1},
\end{equation}
which is associated with  $\mathfrak C(\gamma)$.  By $\mathcal J(\delta)$ we denote a collection of disjoint intervals of length $l\in (2^{-1}\delta, 2\delta)$ which are included in $\mathbb I_0$.  
 For a given function $g$ on $\mathbb I_0\times \mathbb I$ and  $J\in \mathcal J(\delta)$, we set
\[g_J(u,r)=\chi_J(u) g(u,r).
\] 
We denote $\supp_u g= \{ u: \supp f(u, \cdot) \neq \emptyset\}$, so $\supp_u g_J$ is included in  $J$.

Before introducing the decoupling inequality, we define a weight function for a general ball $B\subset\R^m$. Let $c_B$, $R_B$ denote the center of $B$, and the radius of $B$, respectively. Then, $\omega_B(x)=(1+R_B|x-C_B|)^{-N}$ with a sufficiently large $N\geq 100m$. Using the decoupling inequality for the nondegenerate curve  \cite{BDG} and  the argument in \cite{BD} (see also \cite{BGHS2}), we have the following.

\begin{thm}\label{BDGcone} Let  $p \ge d(d+1)$ and $\alpha_d(p):=(2p-d^2-d-2)/(2dp)$. Let $0<\delta<1$ and $\mathcal J(\delta^{1/d})$ be a collection of disjoint intervals given as above.  Let  $B$ denote a ball of radius $\delta^{-1}$ in $\R^{d+1}$.  Suppose that $\gamma$ is nondegenerate. 
Then, for any  $\epsilon>0$  we have
\[\|E^\gamma (\sum_{J\in \mathcal J({\delta^{1/d}})} g_J)\|_{L^p(\omega_{B})}\lesssim_\epsilon \delta^{-\alpha_d(p)-\epsilon}\Big(\sum_{J\in \mathcal J({\delta^{1/d}})}
\|E^\gamma g_J\|_{L^p(\omega_{B})}^p\Big)^{1/p}.\]
\end{thm}

However, the phase functions $\Phi_{\pm}^\theta(x,t,s,\xi)$ is not linear in $t,s,\theta$. So, for our purpose  of proving the smoothing estimate, we need a variable coefficient generalization of  Theorem \ref{BDGcone}.

\subsection{Variable coefficient decoupling} Let
\begin{equation}\label{def of D}
\mathbb D= \mathbb B^{d+1}(0,2)\times (-1,1)\times \mathbb (1/2,2).
\end{equation}
Let  $\Phi :\tilde{\mathbb D} \rightarrow \R$ be a smooth function on  $ \mathbb B^{d+1}(0,2)\times (-1,1)$ and 
$A$ be a smooth function  with $\supp A\subset \mathbb D.$
 For $\lambda\ge1$, we consider 
\[ \mathcal  E_\lambda g(z)=\iint e^{i\lambda r\Phi(z,u)}A(z,u,r)g(u,r)dudr.\]
The following is a variable coefficient generalization of Theorem \ref{BDGcone}. Let  \[\mathcal T(\Phi)(z,u)= (\Phi(z,u), \partial_u\Phi(z,u),\cdots,\partial_u^d\Phi(z,u)).\]

\begin{thm}\label{var decoupling}
Let $p\ge d(d+1)$ and $\mathcal J=\mathcal J(\lambda^{-1/d})$. Suppose that 
\begin{align}
\label{conddec1}
\operatorname{rank}   D_z \mathcal T(\Phi)=d+1
\end{align}
on $\supp a$.  Then, for any $\epsilon>0$  and $M>0$, we have
\Be 
\label{dec}
\|(\sum_{J\in \mathcal J} \mathcal E_\lambda g_J)\|_{L^p}\lesssim_{\epsilon,M} \lambda^{\alpha_d(p)+\epsilon}
\Big(\sum_{J\in \mathcal J} \| \mathcal E_\lambda g_J\|_{L^p}^p\Big)^{1/p}
+\lambda^{-M}\| g\|_2.
\Ee
Here,  we allow discrepancy between amplitude functions  in the left hand and right hand sides, that is to say, the amplitude functions on the both sides are not necessarily the same.  
\end{thm}

It is clear that the model case phase $\Phi(z,u)=z\cdot (1,\gamma(u))$  satisfies the condition \eqref{conddec1} provided that $\gamma$ is nondegenerate.  We will use Theorem \ref{var decoupling}  to prove  Lemma \ref{0->4} and  \ref{4->3}  which   provide the key ingredients for the proof of Proposition \ref{main u3}. By performing allowable  transformations on the  phase $\Phi^\theta_\pm$ (\eqref{th}), it will be verified that  the consequent phase function satisfies  the condition \eqref{conddec1} (see Section \ref{sec-3para}).   Furthermore, note that $\lambda$ in Theorem \ref{var decoupling} plays a role of $\delta^{-1}$ in Theorem \ref{BDGcone}. To see this, one only needs to change variables $z\rightarrow z/\lambda$ in Theorem \ref{var decoupling}.

We refer to the inequality \eqref{dec} as a decoupling of $\E_\lambda$ at scale $\lambda^{-1/d}$.  As is clear, the implicit constant in \eqref{dec} 
is independent of particular choices of $\mathcal{J}(\lambda^{-1/d})$.  The role of the amplitude function $A$ is less significant. 
In fact, changes of variables $z\rightarrow \mathcal Z(z)$ and  $u \rightarrow \mathcal U(u)$ separately in $z$ and $u$  do not have  effect  on the decoupling inequality as long as $\mathcal Z,$ $ \mathcal Z^{-1}$, $\mathcal U$, and $\mathcal U^{-1}$ are smooth with uniformly bounded derivatives up to some large order.  The decoupling for the original operator can be recovered by undoing the changes of variables. This makes it possible  to decouple an operator by using the decoupling inequality in a normalized form.   For our purpose it is enough to consider the amplitude of the form $A(z,u,r)=A_1(z)A_2(u,r)$.  
This can be put together  with those aforementioned changes of variables to obtain the desired decoupling inequalities.  

Theorem \ref{var decoupling}  can be shown through routine adaptation of  the argument in \cite{BHS}, where  the authors obtained a variable coefficient generalization of the decoupling inequality for conic hypersurfaces,  that is to say, Fourier integral operators. However, we include a proof  of Theorem \ref{var decoupling} for convenience of  the readers (see Appendix \ref{proof-var}).

\subsection{Decoupling with a degenerate phase}
To show Proposition \ref{main u2} we also need to consider an operator which does not satisfy the nondegenerate condition  \eqref{conddec1}.  In particular, we make  use of the following for the purpose.

\begin{cor}\label{degendecoup}
Let $\mathcal{J}(\lambda^{-1/d})$ be a collection of disjoint intervals  such that 
$J\subset (-\epsilon_0, \epsilon_0)$ for all $J\in \mathcal{J}(\lambda^{-1/d})$. 
Suppose that
$\det  D_z\mathcal T(\Phi)(z,0)=0$ and 
\Be
\label{degencond}
\det\big( \nabla_{\!z}\Phi(z,0), \partial_u\nabla_{\!z}\Phi(z,0),\cdots,\partial_u^{d-1}\nabla_{\!z}\Phi(z,0),\partial_u^{d+1}\nabla_{\!z}\Phi(z,0)\big)\neq0
\Ee
for $z\in \supp_z A$. Then, if $\epsilon_0$ is small enough, for $\epsilon>0$ and $M>0$ we have
\[
\|(\sum_{J\in\mathcal J(\lambda^{-1/d})} \mathcal  E_\lambda g_J)\|_{L^p}\lesssim_{\epsilon,M} \lambda^{\alpha_d(p)+\epsilon}
\Big(\sum_{J\in \mathcal J(\lambda^{-1/d})}\|\mathcal E_{\lambda}g_J\|_{L^p}^p\Big)^{1/p}
+\lambda^{-M}\|g\|_2.
\]
\end{cor}

We will use Corollary \ref{degendecoup}  to prove  Lemma \ref{2paracase}, which gives the key inequality needed for the proof of Proposition \ref{main u2}. After suitable transformations on $\Phi_{\pm}^0$ (\eqref{th0}), it can be verified that the condition \eqref{degencond} for $d=3$ is satisfied by the resulting phase function (see Section \ref{sec-2para}). Futhermore, Corollary \ref{degendecoup} can be generalized under some types of finite type conditions but we do not intend to pursue this here.
 
 Before beginning the proof of Corollary \ref{degendecoup}, we provide a brief overview of our arugment.   
  A typical example of the phase that satisfies \eqref{degencond} is
 \begin{equation}\label{degenerate model}
     \tilde \Phi_0(z,u):=z\cdot\big(1,u,\cdots,u^{d-1}/(d-1)!,u^{d+1}/(d+1)!\big), 
 \end{equation}
which becomes  nondegenerate away from the origin. This observation can be exploited using dyadic decomposition  and a standard rescaling argument.   Let $j_0$ be the largest integer satisfying 
\[ 2^{j_0}\leq \lambda^{{1}/(d+1)-\epsilon}.\]
We decompose $(-1,1)$ into  dyadic intervals $[-2^{-j_0},2^{-j_0}]$, and $\pm [2^{-j},2^{-j+1}),$ $1\le j\le j_0$.   
For $[-2^{-j_0},2^{-j_0}]$, we just use the trivial decoupling inequality (see \eqref{trivial} below).   
For the intervals $\pm [2^{-j},2^{-j+1}),$ $1\le j\le j_0$, we stretch out $\pm [2^{-j},2^{-j+1})$ to $\pm(1,2]$ by rescaling  so that the phase   $\tilde \Phi_0$ remains unchanged  and 
becomes nondegenerate on $\pm(1,2]$.  Thus, we can apply Theorem \ref{var decoupling}. After undoing  rescaling, 
we apply the  trivial decoupling to get the desired inequality. 
For a general phase, the phase may change after applying the rescaling argument. However, this does not pose an issue, as the decoupling inequality holds uniformly as long as the phase remains sufficiently controlled.

\begin{proof}[Proof of Corollary \ref{degendecoup}]
For $j<j_0$, we set 
\[ g_j=\sum_{  1\le 2^{j} \dist(0, J)<  2 }  g_J .\] 
 Thus, $\sum_J g_J= \sum_{1\le j<j_0}g_j + \sum_{\dist(0,J)<2^{-j_0}}  g_J$.  Let $A_j(z,u,r)=A(z,2^{-j}u,r)$ and $\tilde g_j=2^{-j} g_j(2^{-j}\cdot,\cdot) $. For  
$0\le j<j_0$, changing variables $u\to 2^{-j}u$,  we get
\[
 {\textstyle  \sum_J  \E_{\lambda}g_J=\sum_{0\le j< j_0}}
 \E^{\lambda \Phi(\cdot,2^{-j}\cdot)}_{A_j}  \tilde g_j + {\textstyle \sum_{\dist(0,J)<2^{-j_0}}} \E^{\lambda \Phi}_A g_J.
\]
Here and afterward, for given  $\Psi$ and $b$,  we denote   
\begin{equation}\label{E^psi_b def}
\E^\Psi_b g(z)=\iint e^{ir\Psi(z,u)}b(z,u,r) g(u,r)dudr.
\end{equation}

We set $\tilde \Phi(z,u)=\sum_{k=0}^{d+1}\partial_u^k\Phi(z,0){u^k}/{k!}$. 
Using Taylor's expansion,  we have
\[
\Phi(z,u) = \tilde \Phi(z,u) +\mathcal R(z,u),
\]
where 
$ \mathcal R(z,u)=
\int_0^{u}\partial_u^{d+2}\Phi(z,s)(u-s)^{d+1}ds/{(d+1)!}$. 
From the condition \eqref{degencond} we note that  the vectors $\nabla_{\!z}\Phi(z,0), \partial_u\nabla_{\!z}\Phi(z,0),\cdots,\partial_u^{d-1}\nabla_{\!z}\Phi(z,0)$ are linearly independent. Meanwhile, since $\det  D_z\mathcal T(\Phi)(z,0)=0$, 
$\nabla_{\!z}\Phi(z,0),$ $\cdots,\partial_u^{d-1}\nabla_{\!z}\Phi(z,0),$ $\partial_u^{d}\nabla_z\Phi(z,0)$ are linearly dependent. Thus,  there are smooth functions $r_0, \dots, r_{d-1}$,    such that 
$\partial_u^{d}\Phi(z,0)/d!=  \sum_{k=0}^{d-1} r_k(z)\partial_u^k\Phi(z,0)/k!.$
This yields 
\[ \tilde \Phi(z,u)=   \sum_{k=0}^{d-1}  (  1+ u^{d-k} r_k(z)  )  \partial_u^k\Phi(z,0)\frac{u^{k}}{k!}  +   \partial_u^{d+1}\Phi(z,0)\frac{u^{d+1}}{(d+1)!} .  \]

Let $\mathcal L$ denote  the inverse of  $z \mapsto (\Phi(z,0),\cdots,\partial_u^{d-1}\Phi(z,0),\partial_u^{d+1}\Phi(z,0))$.
Setting   $T_j(z)=\mathcal L (2^{-(d+1)j}z_1,\cdots, 2^{-2j}z_{d}, z_{d+1})$, we have
\[   2^{(d+1)j} \tilde \Phi(T_j (z),2^{-j}u)=  \sum_{k=0}^{d-1}  (  1+ (2^{-j}u)^{d-k} r_k(T_j (z))  )\frac{ z_{k+1}u^{k}}{k!} +\frac{z_{d+1}u^{d+1}}{(d+1)!}.\] 
One can easily see  that  $|\partial^\alpha_{z,u} ( 2^{(d+1)j} \tilde \Phi(T_j (z),2^{-j}u)- \tilde \Phi_0(z,u))|\lesssim 2^{-j}$ and $|\partial^\alpha_{z,u} (2^{(d+1)j} \mathcal R(T_j(z),2^{-j}u))| \lesssim 2^{-j}$ for any multi-index $\alpha$.  
Therefore, 
\[  [\Phi_j](z,u)   :=   2^{(d+1)j}  \Phi(T_j (z),2^{-j} u), \]
which is close to $\tilde \Phi_0(z,u)$ 
defined by \eqref{degenerate model}, 
satisfies the nondegeracy  condition  \eqref{conddec1}  for $|u|\sim 1$   if $\epsilon_0$ is small enough. 
Changing variables $z\to T_j (z)$, we have 
\[  \E^{\lambda \Phi(\cdot,2^{-j}\cdot)}_{A_j}  \tilde g_j (T_j (z))=  \E^{\lambda 2^{-(d+1)j}   [\Phi_j]}_{A_j\circ T_j} \tilde g_j (z).   \]
Decomposing $A_j\circ T_j $ into smooth functions which are supported in a ball of radius $\sim 1$, we apply Theorem  \ref{var decoupling}. By putting together the resultant inequalities on 
each ball, this gives decoupling of $ \E^{\lambda 2^{-(d+1)j}   [\Phi_j]}_{A_j\circ T_j} \tilde g_j $ at scales $\lambda^{-1/d} 2^{(d+1)j/d}$.  Here, it should be noted that 
the constants in the decoupling inequality can be taken uniformly since  the phases  $[\Phi_j]$ are close to $\tilde \Phi_0$. 

 After undoing the change of variables and rescaling it in turn gives decoupling of $\E_\lambda g_j$ at scales $\lambda^{-1/d} 2^{j/d}$.  Now, in order to obtain decoupling at scale $\lambda^{-1/d}$, we make use of the trivial decoupling, that is to say, 
 \Be 
 \label{trivial}
 \|\sum_{J\in\mathcal J} \mathcal E_\lambda g_J\|_{L^p}\le (\#\mathcal J)^{1-2/p}  (\sum_{J\in\mathcal J} \|\mathcal E_\lambda g_J\|_{L^p}^p)^{1/p}
 \Ee  for any collection $\mathcal J$ of disjoint intervals.
Since there are  as many as $\sim 2^{j/d}$ intervals $J$, it produces a factor of $O(2^{j(1-2/p)/d})$ in its bound.
Putting everything  together, we see that $\|\E_\lambda g_j\|_{L^p}$ is now bounded above by a constant times 
\begin{align*}
(\lambda 2^{-(d+1)j})^{\alpha_d(p)+\epsilon} 2^{j(\frac1d-\frac 2{pd})}\Big(\!\sum_{ 1\le   2^{j} \dist(0, J)<  2 }\| \mathcal E_\lambda g_J\|_p^p\Big)^{1/p}  + (2^{(d+1)j}/\lambda)^{M}\| g\|_2.
\end{align*} 
The terms with $\dist(0,J)<2^{-j_0}$ can be handled easily. Since $-{(d+1)\alpha_d(p)}+(1- 2/p)/d<0$, taking summation along $1\le j\le j_0$, 
we get the desired inequality.
\end{proof}

\section{Local smoothing estimates }
\label{nondegeneracy section}
In this section, we prove  Proposition \ref{main u2}  and \ref{main u3} making use of the key observation that 
the immersions \eqref{th0} and \eqref{th} give conic extensions of  finite type curves. 
Using  suitable forms of decoupling inequalities, we first decompose  the  averaging operators so that  the consequent operators have their Fourier supports in narrow angular sectors. For each of those operators, fixing some variables, we make use of   
 the local smoothing estimate for the 2-d wave propagator in $\R^{2+1}$ (for example, see \eqref{st4}), or lower dimensional decoupling inequality.  
 
 Throughout this section, we assume 
\[
\supp \widehat f \subset \mathbb A_j.
\]
To exploit the decoupling inequalities, we decompose $f$ into functions of which Fourier supports are contained in angular sectors.  
For $\kappa\in(0,  1)$,  let   $\{ \Theta^{\kappa}_m \}_{m=1}^N$ denote a collection of disjoint arcs of  length $L\in (2^{-1}\kappa, 2\kappa)$ such that 
$\bigcup_{m=1}^N  \Theta^{\kappa}_m=\mathbb S^1$.  Let $\{\zeta_m^\kappa\}_{j=1}^N$ be a partition of unity on $\mathbb{S}^1$
satisfying $\supp\zeta_m^{\kappa}\subset \Theta_{m-1}^{\kappa}\cup \Theta_{m}^{\kappa}\cup \Theta_{m+1}^{\kappa}$ for $1\le m\le N$ (here, we identify $\Theta_{0}^{\kappa}=\Theta_{N}^{\kappa}$ and 
$\Theta_{N+1}^{\kappa}=\Theta_{1}^{\kappa}$) and $|(d/d\theta)^l\zeta_m^{\kappa}|\lesssim \kappa^{-l}$ for $l\ge 0$. We denote 
\[ \mathfrak{S}(\kappa)=\{  \zeta_m^{\kappa} \}_{j=1}^N.\]
For each $\nu \in \mathfrak{S}(\kappa)$, set 
\[ \widehat{f_{\nu}}(\xi)=\widehat{f}(\xi)\nu (\xi/|\xi|).\]
 
\subsection{2-parameter case: Proof of Proposition \ref{main u2}}
\label{sec-2para}
We only consider the estimate for 
\[  \mathcal{U}^0:= \mathcal{U}^0_+.\]  The estimate for $ \mathcal{U}^0_-$ follows by the same argument. 
We begin with the next lemma, which we obtain by using Corollary \ref{degendecoup}.  
\begin{lem}\label{2paracase}
    Let  $p\geq 12$ and $j\geq 0$. Suppose that $\supp \widehat f\subset \mathbb A_j$.  Then, for any $\epsilon>0$ and $M>0$, we have
    \begin{equation}\label{2para first decomp}
        \Vert \mathcal{U}^0 f\Vert_{L^p}\lesssim_{\epsilon, M} 2^{(\frac{1}{3}-\frac{7}{3p}+\epsilon)j}\Big(\sum_{\nu \in \mathfrak{S}(2^{-j/3})}\Vert \mathcal U^0 f_{\nu} \Vert_{L^p}^p\Big)^{{1}/{p}}+2^{-Mj}\Vert f\Vert_{L^p}.
    \end{equation}
\end{lem}
\begin{proof}  
Decomposing $f$ in the Fourier side, by symmetry  we assume that $\supp \widehat{f}$ is additionally  included in the set  $\lbrace \xi : |\xi_2|\leq 2\xi_1 \rbrace.$
 We make changes of variables  $\xi\to 2^j\xi$ and  $(\xi_1,\xi_2)\to (r,ru)$, successively, to obtain
    \[ \mathcal U^0f(x,t,s, u)= a(x,t,s)\int e^{i2^j r\Phi(x,t,s,u)}\widehat{f(2^{-j}\cdot)}(r, ru) rdrdu, \]
    where $ \Phi(x,t,s, u)=x_1+x_2u+ |(1,u)_{t,s}|$. 
    Let us set 
\Be  
\label{h-def}
h(u)=(\rho^2+u^2)^{-1/2}, \quad   \rho=t/s. 
\Ee
Then, recalling \eqref{comp}, we see that  
 \[ \nabla_{x,t,s}\Phi(x,t,s, u)=\bar\gamma(u):=(1, \gamma(u)):=  \big(1, u,\rho h(u),  {u^2}h(u)\big).   \]
     
 \begin{lem}\label{2p-lem}
 Let $t,s\in \mathbb I$. 
 Then,  we have 
 \begin{align}
\label{33}
 | \det(\bar\gamma(u), \bar\gamma'(u),\bar\gamma''(u), \bar\gamma'''(u)) |\sim |u|,
\\  
\label{44}
 |\det(\bar\gamma(u), \bar\gamma'(u),\bar\gamma''(u),\bar\gamma''''(u))|_{u=0}\sim 1.
\end{align}
\end{lem}

 \begin{proof}  
 Note that
  \[ \gamma^{(k)} (u)=\big (0,\ \rho h^{(k)}(u),\ 2(2k-3)h^{(k-2)}(u)+ 2kuh^{(k-1)}(u)+u^2h^{(k)}(u)\big)\] 
  for $k=2,3$.  Since  $\det(\bar\gamma(u), \bar\gamma'(u),\bar\gamma''(u), \bar\gamma'''(u))=\det(\gamma'(u),\gamma''(u),\gamma'''(u))$, 
   \begin{align*}
   \det(\bar\gamma(u), \bar\gamma'(u),\bar\gamma''(u), \bar\gamma'''(u))
    &=2\rho\det\begin{pmatrix}
      h'' & h+2uh'\\
     h''' & 3h'+3uh''
\end{pmatrix}.
\end{align*}
  After a computation one can easily check  that the following hold: 
\Be
\label{h-deriv}
\begin{aligned}
     h'&(u)=-u h^3(u),
    \qquad \qquad \qquad 
         h''(u)=(2u^2-\rho^2) h^5(u) ,
     \\
         h'''(u)&= 3(3\rho^2u-2u^3)h^7(u), 
       \ \, h''''(u)=3(8u^4-24\rho^2u^2+3\rho^4)h^9(u).
   \end{aligned}
\Ee

Using this, we obtain $\det(\gamma'(u),\gamma''(u), \gamma'''(u))=-{6\rho^5u}(u^2+\rho^2)^{-5}.$  This gives \eqref{33} since $\rho\sim 1$. 
Furthermore, differentiating both sides of the equation, we also have 
$
    \det(\gamma'(u),\gamma''(u),\gamma''''(u))
    =6\rho^5(u^2+\rho^2)^{-6}(9u^2-\rho^2), 
$ 
which shows \eqref{44}. 
\end{proof}

Lemma \ref{2p-lem} shows that \(\nabla_{x,t,s}\Phi(x,t,s,u)\) satisfies the assumption in Corollary \ref{degendecoup} for \(d=3\). Thus, if \(u\) is away from \(u=0\), then \(\nabla_{x,t,s}\Phi(x,t,s,u)\) fulfills the nondegeneracy condition \eqref{conddec1}. 
Therefore, by decomposing the integral \(\mathcal{U}^0 f\) into two parts: one over \(u \in (-\epsilon_0, \epsilon_0)\) for a sufficiently small \(\epsilon_0\) and the other over \(u \in (-\epsilon_0, \epsilon_0)^c\), we can apply Corollary \ref{degendecoup} and Theorem \ref{var decoupling} to the former and the latter, respectively. This allows us to obtain decoupling at scale \(2^{-j/3}\) for both parts.

Note that the $u$-support of $g_\nu(u,r):=\mathcal F{f_\nu(2^{-j}\cdot)}(r,ru)$, $\nu \in \mathfrak{S}(2^{-j/3})$ are contained in boundedly overlapping intervals of length $\sim 2^{-j/3}$, so we have the decoupling inequality 
 \[\|\sum_{\nu \in \mathfrak{S}(2^{-j/3})} \mathcal E_{2^j}  g_\nu\|_p \lesssim 2^{(\frac{1}{3}-\frac{7}{3p}+\epsilon)j} (\sum_{\nu \in \mathfrak{S}(2^{-j/3})}\Vert \mathcal E_{2^j}  g_\nu \Vert_{L^p}^p)^{{1}/{p}} +2^{-Mj}\Vert f\Vert_{L^p}.\]  Therefore,  undoing the changes of variables $\xi\to 2^j\xi$ and  $(\xi_1,\xi_2)\mapsto (r,ru)$, we get the desired inequality  \eqref{2para first decomp}.
\end{proof}
    
To complete the proof of Proposition \ref{main u2} it is sufficient to show \eqref{2-local} for $p> 12$ since the estimate for $4\le p\le 12$ follows by interpolation with the estimate  \eqref{st4}. 
 By the inequality \eqref{2para first decomp}, we only have to prove that
    \[ \Big(\sum_{\nu\in\mathfrak{S}(2^{-j/3})}\Vert \mathcal U^0f_{\nu} \Vert_{L^p}^p\Big)^{1/{p}}\lesssim 2^{(\frac{1}{6}-\frac{2}{3p}+\epsilon)j}\Vert f\Vert_{L^p} \]
    for $p> 12$.  Since $\supp \widehat f\subset \mathbb A_j$, one can easily see that $(\sum_\nu \|f_{\nu} \Vert_{p}^p)^{1/p} \lesssim \|f\|_p$ for $2\le p\le \infty$. As before, this in fact follows by interpolation between the estimates for $p=2$ and 
    $p=\infty$. 
   So, the matter is reduced to showing that 
   \Be 
   \label{hoho} \Vert \mathcal U^0f_{\nu} \Vert_{L^p_{x,t,s}}\lesssim 2^{(\frac{1}{6}-\frac{2}{3p}+\epsilon)j}\Vert f_\nu \Vert_{L^p}\Ee
 for $p\geq 12$.  Recalling \eqref{utheta} with $\theta=0$ and changing variables $\xi_2\to \xi_2/s$, we use the local smoothing estimate for the wave operator. 
    Since $s\sim 1$, the support of $\widehat f_{\nu}(\xi_1,\xi_2/s)$ is included in an angular sector of angle $\sim 2^{-{j}/{3}}$. Applying Lemma \ref{lem:locals} with $\lambda=2^j$ and $b\sim 2^{-{j}/{3}}$, we obtain
    \[ \Vert \mathcal U^0f_{\nu}(x, t, ts) \Vert_{L^p_{x,t}}\le C 2^{(\frac{1}{6}-\frac{2}{3p}+\epsilon)j}\Vert f_\nu \Vert_{L^p}\]
    for any $\epsilon>0$. Integrating in $s$  gives the desired estimate \eqref{hoho}.

\subsection{3-parameter case:  Proof of Proposition \ref{main u3}}
\label{sec-3para}
As before, we only consider the estimates for 
\[  \mathcal{U}^\theta_k:= \mathcal{U}^\theta_{+,k}\] given by \eqref{udef1}.
Those  for $\mathcal{U}^{\theta}_{-,k}$ can be obtained in the same manner. 
We use the decoupling inequality in  Theorem \ref{var decoupling}.  We start with the next lemma.   

\begin{lem}\label{0->4}
    Let $p\geq 20$ and $0\leq k\leq j$, and set $z=(x,t,s,\theta)$. Suppose  $\supp \widehat f\subset \mathbb A_j$. Then, for any $\epsilon>0$ and $M>0$, we have
    \begin{equation}
   \label{e0->4}
    \Vert  \mathcal U^\theta_k f\Vert_{L^p_z}\lesssim_{\epsilon, M} 2^{(1-\frac{11}{p}+\epsilon)\frac{j-k}{4}}\Big(\sum_{\nu\in \mathfrak{S}(2^{({k-j})/{4}})}\Vert \mathcal U^\theta_k  f_{\nu}          
    \Vert_{L^p_z}^p\Big)^{{1}/{p}}+  \mathcal R,
\end{equation}
where $\mathcal R=2^{-M(j-k)}\Vert f\Vert_{p}$. 
 \end{lem}

\begin{proof}
By rotational symmetry, we may assume that  $\theta$ is restricted near $\theta=0$.  Thus, we only need to consider 
    \[    \tilde {\mathcal U} f(z',\theta)=\int e^{i\Psi(z',\theta,\xi)}  \tilde a(z',\theta) \widehat{f}(\xi)d\xi,\quad   z'=(x,t,s),\]
where 
    \[ \Psi(z',\theta,\xi)=x\cdot\xi+|(R_{\theta}^\ast\xi)_{t,s}|, \quad \tilde a(z',\theta)= a(x,t,s)\phi_{<0}(\theta/\epsilon_0)\psi(2^k|t-s|)\]
    for a small $\epsilon_0>0$. Changing variables  $z'\to 2^{-k}z'$ and $\xi\to 2^j \xi$, we have 
    \[    \tilde {\mathcal U} f(z',\theta)= \int e^{i2^{j-k}\Psi(z',\theta,\xi)} \tilde a(2^{-k}z',\theta) \widehat{f(2^{-j}\cdot)}(\xi)d\xi. \]
 
  We decompose   $\tilde a(2^{-k}z',\theta)=\sum_{n} a_{n}(z',\theta)$ such that  $\supp_{z'} a_{n}$ are included in finitely overlapping  balls of radius $1$ and the derivatives of $a_{n}$ are uniformly bounded.   
We are now reduced to obtaining  the decoupling inequality  for  the operator 
\Be \label{sss}   
\mathcal E(\lambda\Psi,  a_{n})g(z',\theta) \coloneqq \int e^{i\lambda \Psi(z',\theta, \xi)} a_{n}(z',\theta) g(\xi) d\xi   \Ee
with $\lambda:=2^{j-k}$ and $ g:= \mathcal F[f(2^{-j}\cdot)]$ whose  support is included in $ \mathbb A_0$. 
  
We now intend to  apply Theorem \ref{var decoupling}. However,  the cutoff $a_{n}$ is no longer supported in a fixed bounded set, so the constants appearing  the decoupling inequality for $\mathcal E(\lambda\Psi,  a_{n})g$
may differ. To guarantee the the constants are uniformly bounded,  one may consider a slightly modified operator. 
Let $z'_0\in \supp_{z'} a_{n}$. Changing variables 
 $z'\to z'+ z'_0$,  we may replace $a_{n}$,  $g$,  $\Psi$  by  $\tilde a_{n}(z',\theta):= a_{n}(z'+z'_0,\theta)$,  $\tilde g(\xi):= e^{i\lambda \Psi(z'_0, 0, \xi)}  g(\xi)$, 
 \[
 \tilde\Psi_n(z',\theta, \xi):=\Psi(z'+z'_0,\theta,\xi)-\Psi(z'_0,0, \xi),
 \]
 respectively. For our purpose it is enough to consider  $\E(\lambda\tilde\Psi_n,  \tilde a_{n})\tilde g$. One can easily check 
 that the derivatives of $\tilde\Psi_n$ and $ \tilde a_{n}$  are uniformly bounded on $B(0,2)\times \mathbb A_0$ for each $n$.

  In order to apply Theorem \ref{var decoupling} to $\mathcal E(\lambda\tilde\Psi_n,  \tilde a_{n})\tilde g$,    we need to verify that  the assumption of Theorem \ref{var decoupling}
 is satisfied after suitable decomposition and allowable transformations.  
Let  $ \mathbb A'= \{(\xi_1,\xi_2)\subset \mathbb A_0:  |\xi_2| < 2\xi_1  \}$ and $\mathbb A''= \{(\xi_1,\xi_2)\subset \mathbb A_0:  |\xi_1| < 2\xi_2  \}.$ Decomposing $g$, we separately consider the following four  cases: 
\Be 
\label{cases}
\supp g \subset  \mathbb A', \ \ \supp g \subset   \mathbb A'', \ \ \supp g \subset - \mathbb A' , \ \  \supp g \subset -\mathbb A''.
\Ee

We first handle the case $\supp g \subset  \mathbb A'$. Writing $\tilde\Psi_n(z',\theta, \xi) = \xi_1 \tilde\Psi_n(z',\theta, 1, \tfrac { \xi_2}{\xi_1})$,   we set 
  \[   \Phi_n(z,u)= \tilde\Psi_n(z, 1, u), \quad z=(z',\theta).\] 
 By recalling \eqref{E^psi_b def}, we see that,  as in the proof of Lemma \ref{2paracase},
 the desired decoupling inequality follows once we obtain a decoupling inequality for 
$\E^{\lambda\Phi_n}_{\tilde a_{n}}$ of scale $\lambda^{-1/4}$.

Even though we have translated $z'\to z'+z_0$, it is more convenient to do computation on $\supp a_n$ before the translation, that is to say, 
$|s-t|\sim 1$ and $t,s \sim 2^k$.  From \eqref{comp} we note  that
\[
\nabla_{z}\Psi(z',\theta, \xi)
:=\Big(\xi,\, \frac{t(R_{\theta}^\ast\xi)_1^2}{|(R_{\theta}^\ast\xi)_{t,s}|},\,\frac{s(R_{\theta}^\ast\xi)_2^2}{|(R_{\theta}^\ast\xi)_{t,s}|},\, 
\frac{2(t^2-s^2)  (R_{\theta}^\ast\xi)_1(R_{\theta}^\ast\xi)_2}{|(R_{\theta}^\ast\xi)_{t,s}|}\Big).
\]
To show that the condition \eqref{conddec1} holds, it is sufficient to consider $\theta=0$ since $\supp_\theta  a\subset (-\epsilon_0, \epsilon_0)$. 
We set  
\Be\label{gam}  \Upsilon(u)  =   \big(u,\rho h(u),u^2h(u),2s^{-1}(t^2-s^2)uh(u)\big)\Ee
for $u\in (-2,2)$. 
Recalling \eqref{h-def},  we see that  $\nabla_z  \Psi(z',0,u)=(1,  \Upsilon(u))$. Thus, to verify  \eqref{conddec1} we have only to show that $\Upsilon$ is nondegenerate, i.e.,   
\[H(u):=\det (\Upsilon'(u), \Upsilon''(u),  \Upsilon'''(u),  \Upsilon''''(u))\neq 0.\]

 To do  this, we note that   ${s H(u)}/{2\rho(t^2-s^2)}$ equals
\[   \det \begin{pmatrix} 
 h'' & 2h+4uh'+u^2h'' &2h'+uh''  \\
 h''' & 6h'+6uh''+u^2h'''& 3h''+uh'''  \\
 h''''& 12h''+8uh'''+u^2h'''' & 4h'''+uh'''' 
\end{pmatrix}
=
\det\begin{pmatrix}
 h'' & 2h & 2h' \\
 h''' & 6h' & 3h'' \\
h'''' & 12h'' & 4h''' 
\end{pmatrix}. \]
Therefore, using \eqref{h-deriv}, one can readily see
\begin{align*}
     H(u)
     =W(u,t,s) \det\begin{pmatrix}
    2u^2-\rho^2 & 1 & -2u 
    \\
  -2u^3+ 3\rho^2u & -u  & 2u^2-\rho^2 
  \\
   8u^4-24\rho^2u^2+3\rho^4 & 4u^2-2\rho^2& -8u^3+ 12\rho^2u 
    \end{pmatrix},
\end{align*}
where $W(u,t,s)=36\rho  {(t^2-s^2)}s^{-1} h^{15}(u)$. A computation shows that the determinant equals $\rho^6$, so  we get 
$H(u)=36\rho^7  {(t^2-s^2)}s^{-1} h^{15}(u)$.   Since $t,s\sim 2^k$ and $|t-s|\sim 1$ on $\supp a_n$, we have
$|H(u)|\ge c$ for a constant $c>0$  on $\supp a_n$.  This shows that $\Phi$ satisfies the nondegeneracy condition \eqref{conddec1} on $\supp \tilde a_{n}\times (-2,2)$ (uniformly for each $n$).

Therefore, by Theorem \ref{var decoupling} with $d=4$ we get decoupling of  $\E^{\lambda\Phi_n}_{\tilde a_{n}}$. In fact,   we get $\ell^p$ decoupling  of $\E(\lambda\Phi, \tilde a_{n})\tilde{g}$  into  $\E(\lambda\Phi, \tilde a_{n})(\tilde{g} \nu (\cdot/|\cdot|))$, $\nu \in \mathfrak{S}(\lambda^{-1/4})$.   
 Putting the inequality for each $n$  together and reversing all changes to recover $\mathcal U^\theta_k f_\nu$, we obtain \eqref{e0->4} when $\supp \widehat f\subset   \mathbb A'. $

For the other cases it is sufficient to show that the  nondegeneracy condition is fulfilled after suitable allowable transformations. 
For the case $\supp g \subset   \mathbb A''$ we write  $\tilde\Psi(z',\theta, \xi) = \xi_2 \tilde\Psi(z',\theta, \xi_1/{\xi_2}, 1)$ and set $  \Phi(z,u)= \tilde\Psi(z, u,1)$.
Then, the matter is reduced to decoupling of the operator $\E(\lambda\Phi, \tilde a_{n})$. Note that 
\[ \nabla_z \Phi(z,u)  =   \big(u,1, u^2\tilde h(u), \tilde\rho \tilde h(u),2t^{-1}(t^2-s^2)u\tilde h(u)\big),\]  
where $\tilde \rho=1/\rho$ and $\tilde h(u)= (\tilde \rho^2+ u^2)^{-1/2}$. Changing coordinates,  we only need to show that the curve  
\[\tilde \Upsilon(u):=\big(u, u^2\tilde h(u), \tilde\rho \tilde h(u),2t^{-1}(t^2-s^2)u\tilde h(u)\big)\] is nondegenerate on $\supp \tilde a_n\times (-2,2)$, i.e., $\det(\tilde \Upsilon', \tilde \Upsilon'', \tilde \Upsilon''', \tilde \Upsilon'''')\neq 0$. This can be easily shown by a similar computation as above. Therefore, $\Phi$ satisfies \eqref{conddec1}.   

The remaining two cases  $\supp g \subset - \mathbb A' , $ $ \supp g \subset -\mathbb A''$ can be handled similarly. So, we omit the details. 
\end{proof}

 To prove Proposition \ref{main u3}, we employ a two-step decoupling argument. In fact, it is not enough for our purpose to use  Lemma \ref{0->4} alone,  as the size of the sectors it provides is too large,  leaving the phase $\Psi(z', \theta, \xi)$ with unexploited curvature properties over those sectors.  More precisely,  the nondegeneracy condition with $d=3$ is satisfied  even after fixing the variable $s$ or $t$. Therefore, we may  apply a lower dimensional decoupling inequality (Lemma \ref{4->3} below). This allows us to fully exploit the curvature property. 

To  this end, we focus on a single operator $f\mapsto\mathcal U^\theta f_{\nu}$ with $\nu\in \mathfrak{S}(2^{({k-j})/{4}})$.

\begin{lem}\label{4->3}
Let $p\geq 12$ and $0\leq k\leq j$.  Let $F= f_{\nu}$ for some $\nu\in\mathfrak{S}(2^{({k-j})/{4}})$. Then, for  any $\epsilon>0$ and $M>0$, we have
\[ \Vert \mathcal U^\theta_k F\Vert_{L^p_z}\lesssim_{\epsilon, M} 2^{(1-\frac{7}{p}+\epsilon)\frac{j-k}{12}}
\Big(\sum_{\nu'\in\mathfrak{S}(2^{(k-j)/{3}})}\Vert \mathcal U^\theta_k  F_{\nu'} \Vert_{L^p_z}^p\Big)^{\frac{1}{p}}+2^{-M(j-k)}\Vert F\Vert_{p}. \]
\end{lem}

\begin{proof} 
We fix $s$ and,  then, apply Theorem \ref{var decoupling} and a rescaling argument, i.e., Lemma \ref{scaling lemma} below with $d=3$, $\lambda=2^{j-k}$, and $\mu=2^{({k-j})/{4}}$. 
Following the same lines of argument as in  the proof of Lemma \ref{0->4}, we need to consider the operator given in \eqref{sss} while fixing $s$, that is to say, 
$\E(\lambda\Psi_s, \tilde a_{n}^s)g$ where $\tilde a_{n}^s(x,t,\theta):= \tilde a_{n}(x,t,s,\theta)$ and 
\[   \tilde\Psi_s(x,t,\theta, \xi):= \tilde\Psi(x,t, s,\theta, \xi).\] 
As before, we may assume $\supp_\theta   \tilde{a}_{n}^s \subset (-\epsilon_0, \epsilon_0) $, and we separately handle the four cases in \eqref{cases}. 
It is enough to consider the first case since the other cases can be handled similarly as in the proof of Lemma \ref{0->4}.   
We consider  
\[ \Phi(x,t,\theta, u):=\tilde\Psi_s(x, t, \theta, 1,u).\]
To show the nondegeneracy condition for $ \Phi$, from \eqref{gam} we only have to show  that the curve 
$ (-2,2)\ni u\mapsto (u,\rho h(u),2s^{-1}(t^2-s^2)uh(u))$  is nondegenerate. 
This is clear  from a similar computation as in the proof of Lemma \ref{0->4}.

As mentioned above, we now apply the rescaling argument, specifically Lemma \ref{scaling lemma} in Appendix \ref{proof-var}with  $\mu=\lambda^{-1/4}$ and $R=\lambda^{1-\delta}$ for a sufficiently small $\delta=\delta(\epsilon)>0$.  
Hence, Theorem \ref{var decoupling} gives $\mathfrak{D}^{\lambda\mu^{d},\epsilon}_{R\mu^{d}}\lesssim 1$.
For the definition of  $\mathfrak{D}^{\lambda, \epsilon}_R$,  we refer the reader to \eqref{decoupling constant def}.  By combining  this  and Lemma \ref{scaling lemma}, we  obtain the desired inequality  via the trivial decoupling inequality.
\end{proof}

We are ready to  complete the proof of  Proposition \ref{main u3}.  As before, we only consider $  \mathcal{U}^\theta_k:= \mathcal{U}^\theta_{+,k}$.    The proof is similar with that of Proposition \ref{main u2} since we now have all the necessary decoupling inequalities. It is sufficient to show \eqref{main u3 estimate} for $p\ge 20$ thanks to the   estimate 
  $  \Vert \mathcal U^\theta_k f\Vert_{L^4_{x,s,t,\theta} } \lesssim 2^{\epsilon j}\Vert f\Vert_{L^4} $ which follows from  \eqref{st4}  by taking integration in $\theta$. 
  Interpolation gives     \eqref{main u3 estimate} for $4\le p <20$. Combining Lemma \ref{0->4} and \ref{4->3} gives 
    \[ \Vert \mathcal{U}_{k}^\theta f\Vert_{L^p_{z} }\lesssim 2^{(\frac{1}{3}-\frac{10}{3p}+\epsilon)(j-k)}(\sum_{\nu\in \mathfrak{S}(2^{(k-j)/3})}\Vert \mathcal U^\theta_k  f_{\nu}\Vert_{L^p_z}^p)^{\frac{1}{p}}+2^{-M(j-k)}\Vert f\Vert_{p} \]
    for $p\ge 20$.     Since $\sum_{\nu\in \mathfrak{S}(2^{(k-j)/3})}  \|f_\nu\|_p^p\lesssim \|f\|_p^p$, it suffices to show 
    \[  \Vert \mathcal U^\theta_k  f_{\nu}\Vert_{L^p_z} \lesssim 2^{(1-\frac{4}{p}+\epsilon)\frac{j+2k}{6}-\frac{k}{p}}\Vert f_\nu\Vert_{L^p}, \quad  \nu\in \mathfrak{S}(2^{(k-j)/3}). \]
     Note that $t,s\sim1$. Changing  variables $s\to ts $ and recalling \eqref{utheta},  we see 
    \[ \Vert \mathcal U^\theta_k  f_{\nu}\Vert_{L^p_z}^p\lesssim  \iint_{|s-1|\lesssim  2^{-k}}   \iint |\tilde{\mathcal{U}}^{\theta,s}_+f_{\nu}(x,t)|^pdx dt   ds   d\theta. \]
    Since $\supp \widehat{f_{\nu}}$ is included in an angular sector of  angle $\sim 2^{{(k-j)}/{3}}$, a similar argument as before and Lemma \ref{lem:locals}  give 
    $\|\tilde{\mathcal U}^{\theta,s}_\pm  f\|_p \lesssim 2^{(1-\frac{4}{p}+\epsilon)\frac{(j+2k)}{6}}\Vert f_{\nu}\Vert_{L^p}.$
     This yields the desired estimate \eqref{main u3 estimate} for $p\ge 20$.

\appendix

\section{Proof of Theorem \ref{var decoupling}}
\label{proof-var}

To prove  Theorem \ref{var decoupling}, we closely follow the argument in \cite{BHS}. Let us denote 
\[ \gamma_\circ(u)=(1, u,u^2/2!, \cdots,u^d/d!).\] 
After suitable decomposition and moderate scaling,  it is enough to consider a class of phase functions which are close to $z\cdot \gamma_\circ(u)$. 
More precisely,  exploiting the assumption \eqref{conddec1},   we can normalize the phase function such that 
\Be 
\label{dcc1} \begin{aligned}
& |\partial_u^{k}\nabla_{z}\Phi- \partial_u^{k}\gamma_\circ|\le \epsilon_0,   \qquad  \qquad  \quad 0\le k\le d,
\\
&\qquad |\partial_u^{k}\partial_z^\beta\Phi|\le \epsilon_0,  \qquad  d+1\le k \le 4N, \  1\le |\beta|\le 4N
\end{aligned}
\Ee
for a small $\epsilon_0>0$ and some large $N$.   Indeed, decomposing the amplitude function $A$,  we assume that 
 \[ \supp A\subset \mathbb B^{d+1}(w,\rho)\times \mathbb B^{1}(v, \rho)\times (1/2,2)\] for some 
$w, v$ and $\rho\ge \lambda^{-1/d}$.  Changing variables $(z,u)\to( z+ w, u+v)$, we replace  
\[\Phi_{w}^{v}(z,u):=\Phi(z+w, u+v)-\Phi(w, u+v),\]  $A_{w}^{v}(z,u):=A(z+w, u+v)$, and
$g_w^v (u,r):=e^{ir\Phi(w,u+v)}g(u+v,r)$  for 
$\Phi$,  $A$, and $g$, respectively. This is harmless because a decoupling inequality for $\E^{\lambda \Phi_{w}^{v}}_{A_{w}^{v}} g_w^v$ gives the corresponding  one for $\E_\lambda g$ as soon as we  undo the procedure. 
 
 Note that 
$  \mathcal T(\Phi_{w}^{v})(0,u)=0$. Thanks to  \eqref{conddec1}, taking $\rho$ to be small enough, 
 we may also assume  by the inverse function theorem  that  the map $ z\mapsto  \mathcal T(\Phi_{w}^{v})(z,u)$  has a smooth local inverse 
 \[ z\mapsto  \mathcal I_{w}^{v}(z,u)\] in a neighborhood of the origin.
 Using Taylor's theorem, we have
\Be
\label{taylor}
 \Phi_{w}^{v}(z,u)=\sum_{k=0}^d\frac{\partial_u^k\Phi_{w}^{v}(z,0)}{k!} {u^k}+ \frac1{d!}\int_{0}^u\partial_u^{d+1}\Phi_{w}^{v}(z,s)(u-s)^d ds.
 \Ee
Setting   $\dil \mu z=(\mu^d z_1, \mu^{d-1} z_2,\cdots,  \mu z_d,  z_{d+1})$ for $\mu>0$,  we have 
\[ 
  [\Phi_{w}^{v}]_\rho(z, u,r):=\rho^{-d} \Phi_{w}^{v}\big( \mathcal I_{w}^{v}(\dil{\rho} z, 0),  \rho u\big)=  z\cdot \gamma_\circ(u) + \mathcal R_w^v (z,u),  
\]
where 
\[ \mathcal R_w^v (z,u)=  \frac{\rho}{d!}\int_{0}^u \partial_u^{d+1}\Phi_{w}^{v}\big( \mathcal I_{w}^{v}(\dil{\rho} z,0),\rho s\big)(u-s)^d ds.\]
Thus, it follows that 
\[
 \E^{\lambda \Phi_w^v}_{\, A_w^v}  g(\mathcal I_{w}^{v}(\dil{\rho} z,0))= \E^{\lambda \rho^d [\Phi_{w}^{v}]_\rho}_{\ [A_{w}^{v}]_\rho} ([g_w^v]_\rho)(z), 
 \]
where  $[A_{w}^{v}]_\rho (z, u,r):=A_{w}^{v}(\mathcal I_{w}^{v}(\dil{\rho} z,0),\rho u,r)$ and  $[g_w^v]_\rho(u,r):=\rho g_w^v(\rho u, r).$ Taking $\rho$ small enough, we have $ |\partial_u^{k}\partial_z^\beta\mathcal R_w^v| \le \epsilon_0$ for $0\le k,  |\beta|\le 4N$ on $\supp\, [A_{w}^{v}]_\rho$.
 Therefore, making additional decomposition of $[A_{w}^{v}]_\rho$  and translation,  we note that 
$\E^{\lambda \Phi_w^v}_{A_w^v}  g(\mathcal I_{w}^{v}(\dil{\rho} z,0))$ can be expressed as a finite sum of the operators 
\[  \E^{\lambda \rho^d \tilde \Phi}_{\tilde A} \tilde g\]
with  $\tilde \Phi$ satisfying \eqref{dcc1} and $\tilde A\in C_c^\infty (\mathbb D)$ where $\mathbb{D}$ is given by \eqref{def of D}. 
Replacing $\lambda \rho^d$ with $\lambda$, we only need to prove the decoupling inequality for the operator of the above form. 
 For the rest of this section we assume that \eqref{dcc1} holds for $\Phi$.

In order to show \eqref{dec}, we make use of Theorem \ref{BDGcone}. For this purpose we set
\[
\Phi_\lambda(z,u)=\lambda\Phi(z/\lambda, u), \quad A_\lambda(z,u,r)=A(z/\lambda, u,r).
\]
For $1\le R\le \lambda$, we denote by $\mathfrak D^{\lambda,\epsilon}_R$ the infimum over all $\mathfrak D$ for which the inequality 
\begin{equation}\label{decoupling constant def}
\|\E^{\Phi_{\lambda}}_{A_\lambda} g\|_{L^p(B)}\le \mathfrak DR^{\alpha_d(p)+\epsilon}\Big(\sum_{J\in\mathcal J(R^{-1/d})}\|\E^{\Phi_\lambda}_{\tilde A_\lambda}g_J\|_{L^p(\omega_{B})}^p\Big)^{\frac 1p}+R^{2d}\big(\frac \lambda R\big)^{-\frac {\epsilon N}{8d}}\| g\|_2
\end{equation}
holds for any ball $B$ of radius $R$ included  in $\B^{d+1}(0,2\lambda)$,  all $\Phi$ satisfying \eqref{dcc1},  and 
$A\in \mathrm C_c^\infty (\mathbb{D})$ with some $\tilde A=\tilde A(A)\in \mathrm C_c^\infty (\mathbb{D})$ satisfying $\|\tilde A\|_{\mathrm C^N} \le \|A\|_{\mathrm C^N}$.

\subsection{Rescaling}
By a rescaling argument, we have the following. 
\begin{lem}\label{scaling lemma}
	Let $R^{-1/d}<\mu<1\le R\le \lambda$.  Let  $B$ be a ball of radius $R$ included  in $\B^{d+1}(0,2\lambda)$.  Suppose $\{J\}\subset  \mathcal J(R^{-1/d})$ and $J\subset  \B^1(v,\mu)$ for some $v\in [-1,1]$.   If $\mu$ is sufficiently small, then
	\begin{align*}
	\| \sum_J \E^{\Phi_{\lambda}} g_J\|_{L^p(\omega_{B})}\!
 \lesssim \mathfrak{D}^{\lambda\mu^{d},\epsilon}_{R\mu^{d}} 
 ( R\mu^{d})^{\alpha_d(p)+\epsilon}
 (\sum_{J}
 \|\E^{{\Phi}_{\lambda}}g_J\|_{L^p(\omega_{B})}^p\!)^{\frac1p}+ 
 \mathcal R
	\end{align*} 
	where $\mathcal R=\mu^2R^{2d}\big(\lambda/ R\big)^{-{\epsilon N}/{8d}}\|g\|_2.$
\end{lem}
We occasionally drop the amplitude functions, which are generically assumed to be admissible.

\begin{proof}
To prove Lemma \ref{scaling lemma}, we only need to consider $\| \sum_J \E^{\Phi_{\lambda}} g_J\|_{L^p(B)}$ instead of $\| \sum_J \E^{\Phi_{\lambda}} g_J\|_{L^p(\omega_{B})}$.
Since $\omega_{B}$ is bounded by a rapidly decreasing sum of characteristic functions, the bounds on $\| \sum_J \E^{\Phi_{\lambda}} g_J\|_{L^p(B)}$ imply those for $\| \sum_J \E^{\Phi_{\lambda}} g_J\|_{L^p(\omega_{B})}$.

Let $B=\mathbb B^{d+1}(\lambda w, R)$ for some $w$.  We make a slightly different form of scaling from the previous one to ensure that the consequent phase satisfies \eqref{dcc1}.  Recalling \eqref{taylor}, we have 
\[ \lambda \Phi_{w}^{v}( \mathcal I^v_{w} ( D_\mu'  \tfrac z \lambda ,0), \mu u)=
z\cdot \gamma_\circ(u)+ \frac{\mu^{d+1}}{d!}\int_{0}^{u}\partial_u^{d+1}\Phi_{w}^{v}
\big(\mathcal{I}^v_{w}(D_\mu'  \tfrac z \lambda,0),\mu s\big)(u-s)^d ds,\]
where
$D'_{\mu}z=(z_1,\mu^{-1}z_2,\cdots,\mu^{-d}z_{d+1})$. Setting 
\[ \tilde \Phi(z,u,r)=z\cdot \gamma_\circ(u)+  \frac{\mu}{d!}\int_{0}^{u}\partial_u^{d+1}\Phi_{w}^{v}(D_\mu  z,\mu s)(u-s)^d ds,\] 
we have  $ \lambda \Phi_{w}^{v}( \mathcal I^v_{w} ( \mathrm D_\mu'  \tfrac z \lambda,0), \mu u)=\tilde \Phi_{\lambda\mu^d}(z,u)$. 
Thus,  it follows that 
\[ \E^{  \Phi_{w}^{v}}_{\,\, A_{w}^{v}}  g(\lambda\mathcal I_{w}^{v}(\mathrm D_{\mu}' z/\lambda,0))= \E^{\tilde \Phi_{\lambda\mu^d}}_{\tilde A} ([g_{w}^v]_\mu)(z/\lambda)\]
where 
$\tilde {A}(z, u,r)=A_{w}^{v}(\mathcal I_w^v(D_{\mu}'z/\lambda) ,0),\mu u,r)$. Changing variables $(z,u)\to (z+w, \mu u+v)$ and $z\to  \mathcal I(z):= \lambda\mathcal I_{w}^{v}(\dil{\mu}' z/\lambda,0)$ gives 
\Be
\label{rseq1}
\|\sum_J \E^{\Phi_{\lambda}}_{\, A}g_J\|_{L^p(B)}\lesssim 
\mu^{-\frac{d^2+d}{2p}}\|\sum_J \E^{\tilde \Phi_{\lambda\mu^d}}_{\,\tilde A} \tilde g_J\|_{L^p(\mathcal I^{-1}(B))},
\Ee
 where $ \tilde g_J= [(g_J)_{w}^v]_\mu.$  We cover $\mathcal I^{-1}(B)$ by a  collection  $\mathbf B$ of finitely overlapping $R\mu^{d}$-balls. So, we have
\[ \|\sum_J \E^{\tilde \Phi_{\lambda\mu^d}}_{\,\tilde A} \tilde g_J\|_{L^p(\mathcal I^{-1}(B))}  \lesssim \big(\sum_{ B'\in\mathbf B} \| \sum_J \E^{\tilde \Phi _{\lambda\mu^{d} }}_{\tilde A} \tilde g_J\|_{L^p(B')}^p\big)^{1/p}.\]
Here, we note that $\supp_z\tilde a$ may be not included in  $B(0,\mu^d R)$. However, by harmless translation 
we may assume  that $\supp_z\tilde A\subset B(0,\mu^d R)$  by replacing  the phase and amplitude functions with  $[\tilde \Phi]_{w'}^{0}$ and $[\tilde A]_{w'}^{0}$  for some $w'$ since undoing the translation recovers the desired decoupling inequality.

   Note that $\tilde \Phi$ satisfies \eqref{dcc1} if $\mu$ is small enough. Since $\supp_u \tilde g_J$ are included in disjoint intervals of length $\sim\mu^{-1} R^{-1/d}$,  we now have
\begin{align*}
  \| \sum_J \E^{\tilde \Phi _{\lambda\mu^{d} }}_{\tilde A} \tilde g_J\|_{L^p(B')}
\le&
\mathfrak{D}^{\lambda\mu^{d},\epsilon}_{R\mu^{d}}
(R\mu^{d})^{\alpha_d(p)+\epsilon}\Big(\sum_{J}\| \E^{\tilde \Phi _{\lambda\mu^{d} }}_{\tilde A}\tilde g_J\|_{L^p(\omega_{B'})}^p\Big)^{1/p}
+\mathcal R,  
\end{align*}
where $\mathcal R= (R\mu^{d})^{2d}(\lambda/ R)^{- {\epsilon N}/{8d}}\|\sum_J \tilde g_J\|_2$. We put together the inequalities  over each $B'$  and then reverse the various changes of variables so far to recover the operator $\E^{\Phi_{\lambda}}$, possibly with  a different amplitude function. However, 
the decoupling state is not changed.    
Since  $\#\mathbf B\lesssim \mu^{-d(d+1)/2}$,  we can conclude that
$\|\E^{\Phi_{\lambda}}g\|_{L^p(\omega_{B})}$
is bounded by a constant times
\begin{align*}
\mathfrak{D}^{\lambda\mu^{d},\epsilon}_{R\mu^{d}}
(R\mu^{d})^{\alpha_d(p)+\epsilon}&
(\sum_{J}\|\E^{\Phi_{\lambda}}g_J\|_{L^p(\omega_{B})}^p)^{\frac1p}
+\mu^{-\frac{d^2+d}p+\frac12} (R\mu^{d})^{2d}\big(\tfrac {\lambda \mu^d} R  \big)^{-\frac {\epsilon N}{8d}}\|g \|_p
\end{align*}
where  $g=  \sum_J g_J$. 
Finally, using  $(d^2+d)/p\le 2d^2-2$, we can get the desired result. 
\end{proof}

\subsection{Linearization of the phase}
Let $\Phi$ be a smooth phase satisfying \eqref{dcc1}. For simplicity, denote $\partial_k=\partial_{z_k}$, $k=1, \dots, d+1$.  From \eqref{dcc1}, we have 
$\partial_u(\partial_{2}\Phi/\partial_{1}\Phi)-1=O(\epsilon_0)$.
Thus,  there exists the map $\eta_{z}$ 
such that $(\partial_{2}\Phi/\partial_{{1}}\Phi)(z,\eta_{z}(u))=u$.  
Let 
\[
\Gamma_{z}(u)=\frac{\nabla_{\!z}\Phi(z,\eta_{z}(u))}{\partial_{1}\Phi(z,\eta_{z}(u))}. 
\]
 Also note from  \eqref{dcc1} that $\partial_1 \Phi-1=O(\epsilon_0)$. Furthermore, we have $\Gamma_{z}\cdot e_1=1$ and $\Gamma_{z}\cdot e_{2}=u$ 
where  $\{e_1, \dots, e_{d+1}\}$ is the standard basis in $\R^{d+1}$. 

Let $\lambda w \in \B^{d+1}(0,2\lambda)$.  A Taylor expansion and changing variables $u\to \eta_{w}^{-1}(u)$  give
\[
 \Phi_\lambda(z+\lambda w,u)- \Phi_\lambda(\lambda w,u)
=\partial_{{1}}\Phi(w, u) \Gamma_{w}(\eta_{w}^{-1}(u))\cdot z+ \mathcal R^\lambda_{w}(z,u),
\]
where
\[\mathcal R^\lambda_{w}(z, u)=\frac 1 \lambda\int_0^1(1-\tau)\big\langle \operatorname{Hess}_z \Phi\big(\lambda^{-1}\tau z+w, \eta_{w}^{-1}(u)\big) z, z\big\rangle d\tau. \]

Let us set 
\[ \Omega_{z}(u,r)= (\eta_{z}(u), r/\partial_{1}\Phi(z,\eta_{z}(u))).\] Then, \eqref{dcc1} ensures that $\Omega_{z}$ is smooth. 
Changing variables $(u,r)\to  \Omega_{w}(u,r)$, we see that  
$\E^{\Phi_\lambda}_{A_\lambda} g(z+\lambda w)$ is equal to 
 \Be  \label{eee}
\iint e^{ir z\cdot\Gamma_{w}(u)}    A_{w}(z,\Omega_{w}  (u,r))     (g_{w}\!\circ\, \Omega_{w} ) (u,r)   \frac{ \eta_{w}'(u) dudr   }{ \partial_{1}\Phi(w,\eta_{w}(u))} ,  
\Ee
where 
\begin{align*}
A_{w}(z,u,r) =e^{i r\mathcal R^\lambda_{w}(z, u )}\, A(\lambda^{-1}z+w, u,r), \quad  
g_{z}(u,r)=e^{ir \lambda\Phi(z,u)} g(u,r). 
\end{align*} 
 For this  operator   we could directly apply Theorem \ref{BDGcone} if it were not for   the extra factor $e^{i r\mathcal R^\lambda_{w}(z,\eta_{w} (u) )}$.  
 This is not generally allowed. However, if $|z|\le \lambda^{1/2}$,  expanding it into Fourier series  in $(u,r)$, 
 we may disregard  it as  an minor error. 
 
 More precisely, from \eqref{dcc1}  we note that  $\partial_{1}\Phi-1=O(\epsilon_0)$ and  $\eta'_{z}-1=O(\epsilon_0)$. With a sufficiently small $\epsilon_0$ we may assume  that 
$ g_{w}\!\!\circ \Omega_{w}$ is supported in  $(-1,1)\times [1,2]$. 
 Using \eqref{dcc1}, we have $
|\partial_u^k\mathcal R^\lambda_{w}(z,u)|\le C  {|z|^2}/\lambda$ for $0\le k\le 4N.
$
Consequently, 
\Be 
\label{errore}
|\partial_u^k\big(A_{w}(z,\Omega_{w}  (u/r,r)) |\le C, \quad  0\le k\le 4N 
\Ee 
if  $|z|\le \lambda^{1/2}$. 
Thus, expanding $A_{w}(z,\Omega_{w}  (u/r,r))$ into Fourier series, we have 
$A_{w}(z,\Omega_{w}  (u,r))
=\sum_{\ell \in \Z^2} b_\ell(z)e^{ir\ell\cdot(1,u)}$
with 
$|b_\ell (z)|\lesssim_N (1+|l|)^{-N}.$  From \eqref{eee} we have 
\Be
\label{expand}
|\E^{\Phi_\lambda}_{A_\lambda} g(z+\lambda w)|\le \sum_{\ell\in \Z^2}(1+|\ell|)^{-N}|E^{\Gamma_{w}} (\tilde g_{w})(z+ v_\ell )|,
\Ee
for $|z|\lesssim \lambda^{1/2}$ 
where  $v_\ell:= \ell_1e_1+ {\ell_2}e_{2}$ and 
 \[  \tilde g_{w}(u,r)=   { (g_{w}\!\circ\, \Omega_{w})  (u,r)  \eta_{w}'(u)   }/{ \partial_{1}\Phi(w,\eta_{w}(u))}.  \]
This almost allows us to obtain the first part of the next lemma, which is basically the same as Lemma 2.6 in \cite{BHS}. We recall \eqref{def of E}.

\begin{lem}
Let $0<\delta\le 1/2$ and $1\le \rho\le \lambda^{1/2-\delta}$. Let  $B:=\B^{d+1}(\lambda w,\rho)\subset \B^{d+1}(0,3\lambda/4)$ and $B_0:=\B^{d+1}(0,\rho)$. Suppose that $\Phi$ satisfies \eqref{dcc1}. Then 
\Be
\label{approx1}
\|\E^{\Phi_\lambda}_{A_\lambda}g\|_{L^p(\omega_{B})}\lesssim \|E^{\Gamma_{w}} (\tilde g_{w})\|_{L^p(\omega_{B_0})}+\lambda^{-\delta N/2}\| g\|_2.
\Ee
Additionally, assume that $|w|\le \lambda^{1-\delta'}$. Then, for some admissible  $\tilde A$, we have
\Be
\label{approx2}
\|E^{\Gamma_{w}} (\tilde g_{w})\|_{L^p(\omega_{B_0})}\lesssim \|\E^{\Phi_\lambda}_{\tilde A_\lambda}g\|_{L^p(\omega_{B})}+\lambda^{-\min\{\delta,\delta'\} N/2}\| g\|_2,
\Ee
\end{lem}

\begin{proof} 
For \eqref{approx1}, we separately consider two cases $|z-\lambda w|\le\lambda^{1/2}$ and $|z-\lambda w|> \lambda^{1/2}$. 
We first consider the case $|z-\lambda w|> \lambda^{1/2}$. So, we have 
$
\omega_{B}(z)
\lesssim \lambda^{-\delta(N-d-2)}(1+\rho^{-1}|z-\lambda w|)^{-d-2}.
$
Combining this and a trivial inequality $|\E^{\Phi_\lambda}_{A_\lambda} g|\lesssim \| g\|_2$, we have
\Be
\label{errorex}
\| \chi_{B(\lambda w,\lambda^{1/2})^c}\E^{\Phi_\lambda}_{A_\lambda}g\|_{L^p(\omega_{B})}
\lesssim \lambda^{-\delta N/2}\| g\|_2
\Ee
for a sufficiently large $N$.  Next, we handle the remaining  part 
$\chi_{B(\lambda w,2\lambda^{1/2})}\E^{\Phi_\lambda}_{A_\lambda} g$. Using  \eqref{expand} and H\"older's inequality in $l$, one can obtain
\begin{align*}
\|\chi_{B(\lambda w,2\lambda^{1/2})}\E^{\Phi_\lambda}_{A_\lambda} g\|_{L^p(\omega_{B})}
\lesssim 
\Big\|E^{\Gamma_{w}} (\tilde g_{w}) \sum_{\ell} \frac{ \omega_{B(v_l, \rho)}^{1/p}}{(1+|\ell|)^{N}}\Big \|_{L^p} 
\lesssim     \|E^{\Gamma_{w}} (\tilde g_{w}) \|_{L^p(\omega_{B_0})}. 
\end{align*}
 The second inequality follows from the fact that $\sum_{\ell\in \Z^2}(1+|\ell|)^{-N}\omega_{B(v_\ell, \rho)}^{1/p}\lesssim\omega_{B(0, \rho)}^{1/p}$. By the above inequality and \eqref{errorex}, we conclude that \eqref{approx1} holds.

To show \eqref{approx2}, we use a similar argument. By the same reason as in the proof of \eqref{approx1}, we have $
\|E^{\Gamma_{w}} (\tilde g_{w}) (1-\chi_{B(0,2\lambda^{1/2})})\|_{L^p(\omega_{B_0})}
\lesssim \lambda^{-\delta N/2}\| g\|_2 
$
for a sufficiently large $N$. For the integral over the set $B(0,2\lambda^{1/2})$,  we now undo the changes of variables including $(u,r)\mapsto  \Omega_{w}(u,r)$ which are performed to get \eqref{eee}. 
Consequently, we have
\[
E^{\Gamma_{w}} (\tilde g_{w})(z) =\iint e^{ir \Phi_\lambda (z,u)}    \tilde A_{w}(z,u,r)    g(u,r) dudr , 
\]
where $ \tilde A_{w}(z,u,r)=e^{-i r\mathcal R^\lambda_{w}(z, u )}\, a(\lambda^{-1}z+w, u,r)$. As before,  we can expand 
the function $ \tilde A_{w}(z, \Omega_{w}( u/r,r))$ (cf. \eqref{errore}) into Fourier series in $u,r$.  
Thus,  if  $z\in B(0,2\lambda^{1/2})$, 
\[
\textstyle |E^{\Gamma_{w}} (\tilde g_{w})(z) |
\le C_N \sum_{\ell\in \Z^2}(1+|\ell|)^{-2N}
|\E^{\Phi_\lambda}_{\tilde A_\lambda}(g_\ell)(z)|,
\]
for a suitable symbol $\tilde A$ where 
\[ g_\ell:=e^{i \ell\cdot  \tilde \Omega_{w}^{-1}( u,r)}g.\] 

We  again perform the previous linearization procedure for $\E^{\Phi_\lambda}_{\tilde A_\lambda}(g_\ell)$. 
Since $ \tilde \Omega_{w}^{-1}\circ\Omega_{w}(u,r)=(ru,r),$ by \eqref{approx1} we have 
\[
\|\E^{\Phi_\lambda}_{\tilde A_\lambda}(g_\ell) \|_{L^p(\omega_{B})}
\lesssim \|E^{\Gamma_{w}} (\tilde g_{w})\|_{L^p(\omega_{B(v_\ell, \rho)})}+\lambda^{-\delta N/2}\|g\|_2.
\] 
By this inequality we have
$
S:=\sum_{|\ell|\ge M}(1+|\ell|)^{-2N}
\|\E^{\Phi_\lambda}_{\tilde A_\lambda}(g_\ell) \|_{L^p(\omega_{B})}$ bounded by a constant times $
  M^{-N}\|E^{\Gamma_{w}} (\tilde g_{w})\|_{L^p(\omega_{B_0})}
+\lambda^{-\delta N/2}\|g\|_2. 
$
If we choose a sufficiently large $M$,  the part $S$ can be absorbed in the left hand side of \eqref{approx2}.  
Thus, we obtain
\begin{align*}
\textstyle \|E^{\Gamma_{w}} \tilde g_{w}\|_{L^p(\omega_{B_0})}
\lesssim \sum_{|\ell|<M}
\|\E^{\Phi_\lambda}_{\tilde A_\lambda}g_\ell \|_{L^p(\omega_{B(w, \rho)})}
+\lambda^{-\delta N/2}\|g\|_2.
\end{align*}

We note  that  $\E^{\Phi_\lambda}_{\tilde A_\lambda}g_\ell= \E^{\Phi_\lambda}_{\tilde A_{\lambda,\ell}}g$
where  $\tilde A_{\lambda, \ell}=\tilde A_\lambda e^{i \ell\cdot  \tilde \Omega_{w}^{-1}( u,r)}$. 
Expanding $\tilde A_{\lambda, \ell}$ in  a Taylor series one can get  amplitude functions which are independent  of a particular $B$. 
From those one can find an operator which has the  desired property by pigeonholing.  See \cite{BHS} for details. 
\end{proof}

\subsection{Proof of Theorem \ref{var decoupling}}

Assume that $\Phi$ satisfies \eqref{dcc1} and 
\[1\le K\le R \le \lambda^{1-\epsilon/d}.\]  Let 
$\mathcal{J}:=\mathcal J(R^{-1/d})$ be a collection of disjoint intervals. For simplicity we set   $ g=\sum_{J\in \mathcal J(R^{-1/d})} g_J$. Partition $\mathcal J(R^{-1/d})$ in such a way that 
there is a collection  $\mathcal J'$ of disjoint intervals $J'$ of length $\sim K^{-1/d}$ which include
each interval in $\mathcal J(R^{-1/d})$. Consequently, we have 
\[ \textstyle   g= \sum_{J'\in \mathcal J'} g_{J'}=  \sum_{J'\in \mathcal J'} \sum_{J\in \mathcal J: J\subset J'} g_{J}.  \] 

We consider a ball $B$ of radius $R$ included in $B(0,\lambda)$ and  a collection $\mathcal B_K$ of  finitely overlapping balls $B'$ of radius $K=\lambda^{1/4}$ which covers $B$. Since $R \le \lambda^{1-\epsilon/d}$, one may assume that the center of $B'\in \mathcal B_K$ lies in $B(0,\lambda^{1-\epsilon/d})$ after a translation.
Using \eqref{approx1}, we have
\[
 \|\mathcal E_{\lambda} g \|_{L^p(B)}\lesssim (\sum_{B'\in \mathcal B_K} \|E^{\Gamma_{c_{B', \lambda}}} \tilde g_{c_{B',\lambda}}\|_{L^p(\omega_{B(0, K)})}^p)^{\frac 1 p}+\big(\tfrac R K\big)^{d+1}\lambda^{- N/8}\|g\|_2.
\]
Here $c_{B', \lambda}=\lambda^{-1} c_{B'}$ and $c_{B'}$  denotes the center of $B'$.  
We apply Theorem \ref{BDGcone}  to each $B'\in \mathcal B_K$ and  \eqref{approx2} subsequently  to get the decoupling inequality  at scale $K^{-1/d}$. Consequently, combining the inequality on each $B'$, we obtain
\[
\textstyle \|\E_{\lambda} g\|_{L^p(B_R)}\lesssim K^{{\alpha_d(p)+\epsilon}}\big(\sum_{J'\in \mathcal J' }\|\E_{\lambda}g_{J'}\|_{L^p(\omega_{B_R})}^p\big)^{\frac 1p}
+K^{-1}R^{2d}\big(\tfrac \lambda R\big)^{-\frac{\epsilon N}{8d}}\|g\|_2.
\]
Using Lemma  \ref{scaling lemma}, we get
\begin{align*}
\textstyle \|\mathcal E_{\lambda}g\|_{L^p(\omega_{B})}\!\lesssim
&  \mathfrak D^{\lambda K^{-1},\epsilon}_{R K^{-1}}R^{\alpha_d(p)+\epsilon}
\big(\sum_{J\in\mathcal J}\|\E_{\lambda}g_J\|_{L^p(\omega_{B})}^p\big)^{\frac 1 p} +K^{-\frac \epsilon d}R^{2d}(\tfrac \lambda R)^{-\frac{ \epsilon N}{8d}}\|g\|_2. 
\end{align*}
Thus, for a sufficiently large $\lambda$, we have $
\mathfrak D^{\lambda,\epsilon}_{R}
\le \mathfrak D^{\lambda^{ 3/ 4}, \epsilon}_{R \lambda^{-1/ 4}}.
$
Iteratively applying this inequality,  one can show $\mathfrak D^{\lambda,\epsilon}_{\lambda^{1-\epsilon} }\lesssim \lambda^\delta$ for any $\delta>0$, which completes the proof of Theorem \ref{var decoupling}.

\subsection*{Acknowledgement}  The authors would like to thank the referee for careful reading and valuable suggestions.  This work was supported by the National Research Foundation (Republic of
Korea) grant no. 2022R1A4A1018904 (J. Lee and S. Lee), the KIAS individual grant SP089101(S. Oh) and MG098901(J. Lee).

\bibliographystyle{plain}

\end{document}